\documentclass[11pt]{article}
\usepackage{amsthm}
\usepackage{pxfonts}
\usepackage{yfonts}
\usepackage{dsfont}
\usepackage{graphicx}
\usepackage{relsize}
\usepackage{color}

\def\bel{\begin{equation}\label}
\def\eeq{\end{equation}}
\def\ds{\displaystyle}

\def\mt{\longrightarrow}
\def\v{\vskip 1em}

\def\ve{\varepsilon}

\def\Rec{{\bf R}}
\def\R{\mathbb R}
\def\Z{\mathbb Z}
\def\C{\mathfrak{B}}

\def\N{{\bf N}}

\def\S{{\bf S}}

\def\F{{\bf F}}

\def\Q{{\bf Q}}

\def\A{{\bf A}}
\def\B{{\bf B}}

\def\L{{\bf L}}

\def\V{{\bf V}}

\def\T{{\bf T}}

\def\p{{\partial}}
\def\a{{\bf a}}

\def\bar{\overline}

\def\vol{{\bf vol}}

\def\I{{\bf I}}
\def\II{{\bf II}}

\def\alpha{\alphaup}
\def\beta{\betaup}
\def\gamma{\gammaup}
\def\delta{\deltaup}
\def\theta{\vartheta}
\def\xi{{\xiup}}
\def\eta{{\etaup}}
\def\tau{{\tauup}}
\def\rho{{\rhoup}}
\def\phi{{\phiup}}
\def\psi{{\psiup}}
\def\lambda{{\lambdaup}}
\def\omega{\omegaup}
\def\varphi{{\varphiup}}
\def\gamma{{\gammaup}}

\newtheorem{remark}{Remark}[section]

\begin{document}

\[\hbox{\LARGE{\bf Stein-Weiss inequality revisit on Heisenberg group}}\]

\[\hbox{Chuhan Sun~~~~~and~~~~~Zipeng Wang}\]
\begin{abstract}
We study a family of fractional integral operators  defined   as
\[\I_{\alpha\beta\vartheta}f(u,v,t)~=~\iiint_{\R^{2n+1}} f(\xi,\eta,\tau)\V^{\alpha\beta\vartheta}\big[(u,v,t)\odot(\xi,\eta,\tau)^{-1}\big]d\xi d\eta d\tau\]
where $\odot$ denotes the multiplication law of a Heisenberg group.
$\V^{\alpha\beta\vartheta}$ is a distribution  satisfying Zygmund dilation. 
\v
Let $\omega(u,v)={\sqrt{|u|^2+|v|^2}}^{-\gamma}$, $\sigma(u,v)={\sqrt{|u|^2+|v|^2}}^\delta$.
A characterization is established between  $\omega\I_{\alpha\beta\vartheta}\sigma^{-1}\colon\L^p(\R^{2n+1})\mt\L^q(\R^{2n+1})$ and the necessary constraints consisting of $\alpha,\beta,\vartheta, \gamma,\delta\in\R$ for $1<p< q<\infty$. 
\end{abstract}

\section{Introduction}
\setcounter{equation}{0}
To begin, we recall two classical results of fractional integration on  Euclidean space.
Define
\bel{I_a}
\T_\a f(x)~=~\int_{\R^\N}f(y)\left[{1\over|x-y|}\right]^{\N-\a} dy,\qquad 0<\a<\N.
\eeq
In 1928, Hardy and Littlewood \cite{Hardy-Littlewood} first established an regularity theorem for $\T_\a$ for $\N=1$. Ten years later, Sobolev \cite{Sobolev} made extensions on every higher dimensional space.

{\bf Hardy-Littlewood-Sobolev theorem}~~{\it Let $\T_\a$ defined in (\ref{I_a}) for $0<\a<\N$. We have
\bel{HLS Ineq}
\begin{array}{cc}\ds
\left\| \T_\a f\right\|_{\L^q(\R^\N)}~\leq~\C_{p~q}~\left\| f\right\|_{\L^p(\R^\N)},\qquad 1<p<q<\infty
\\\\ \ds
\hbox{\small{if and only if}}\qquad {\a\over \N}~=~{1\over p}-{1\over q}.~~~~~~~
\end{array}
\eeq
}\\
In 1958, Stein and Weiss \cite{Stein-Weiss} obtained a weighted analogue of the above regularity theorem  by considering the {\it weights} to be suitable powers. 

{\bf Stein-Weiss theorem}~~{\it Let $\T_\a$ defined in (\ref{I_a}) for $0<\a<\N$ and $\omega(x)=|x|^{-\gamma},\sigma(x)=|x|^\delta$ for $\gamma,\delta\in\R$ whenever $x\neq0$. We have
\bel{Stein-Weiss Ineq}
\left\|\omega \T_\a f\right\|_{\L^q(\R^\N)}~\leq~\C_{p~q~\gamma~\delta}~\left\| f\sigma\right\|_{\L^p(\R^\N)},\qquad 1<p\leq q<\infty
\eeq
if and only if
\bel{formula}
 \gamma<{\N\over q},\qquad \delta<\N \left({p-1\over p}\right),\qquad \gamma+\delta\ge0,\qquad {\a\over \N}~=~{1\over p}-{1\over q}+{\gamma+\delta\over\N}.
\eeq
}

$\diamond$ {\small Throughout, $\C>0$ is  a generic constant depending on its sub-indices.}

\begin{remark} In the original paper of Stein and Weiss \cite{Stein-Weiss}, (\ref{formula}) is given as a sufficient condition. Conversely, it turns out to be necessary as well. See the appendix in \cite{Sun-Wang}.
\end{remark}

{\bf Hardy-Littlewood-Sobolev theorem} was first re-investigated by Folland and Stein \cite{Folland-Stein} on Heisenberg group.
We work on its real variable representation with a multiplication law:
\bel{multiplication law}
(u,v,t)\odot(\xi,\eta,\tau)~=~\Big[u+\xi, v+\eta,t+\tau+\mu\big(u\cdot\eta-v\cdot\xi\big)\Big],\qquad\mu\in\R
\eeq
for every $(u,v,t)\in\R^n\times\R^n\times\R$ and $(\xi,\eta,\tau)^{-1}=(-\xi,-\eta,-\tau)\in\R^n\times\R^n\times\R$.

Let $0<\rho<n+1$. Consider
\bel{S_a}
\S_\rho f(u,v,t)~=~\iiint_{\R^{2n+1}}f(\xi,\eta,\tau)\Omega^\rho\Big[(u,v,t)\odot(\xi,\eta,\tau)^{-1}\Big]d\xi d\eta d\tau.
\eeq
$\Omega^\rho$ is a distribution in $\R^{2n+1}$ agree with
\bel{Omega^a}
\Omega^\rho(u,v,t)~=~\left[{1\over |u|^2+|v|^2+|t|}\right]^{n+1-\rho},\qquad \hbox{\small{$(u,v,t)\neq(0,0,0)$}}.
\eeq

{\bf Folland-Stein theorem} ~~{\it Let $\S_\rho$ defined in (\ref{S_a})-(\ref{Omega^a}) for $0<\rho<n+1$. We have
\bel{Folland-Stein theorem}
\begin{array}{cc}\ds
\left\| \S_\rho f\right\|_{\L^q(\R^{n+1})}~\leq~\C_{p~q}~\left\| f\right\|_{\L^p(\R^{2n+1})},\qquad 1<p<q<\infty
\\\\ \ds
\hbox{if and only if}\qquad {\rho\over n+1}~=~{1\over p}-{1\over q}.
\end{array}
\eeq}\\
The best constant for the $\L^p\mt\L^q$-norm inequality in (\ref{Folland-Stein theorem}) is found by Frank and Lieb \cite{Frank-Lieb}. A discrete analogue of this result has been obtained by Pierce \cite{Pierce}. Recently, the regarding commutator estimates are established  by Fanelli and~Roncal \cite{Luca-Luz}.

{\bf Stein-Weiss theorem} has been re-investigated on Heisenberg group by Han, Lu and Zhu \cite{Han-Lu-Zhu}.
\v
{\bf Han-Lu-Zhu theorem} ~~{\it Let $\S_\rho$ defined in (\ref{S_a})-(\ref{Omega^a}) for $0<\rho<n+1$. Suppose $\gamma,\delta\in\R$ and $\omega(u,v)={\sqrt{|u|^2+|v|^2}}^{-\gamma}, \sigma(u,v)={\sqrt{|u|^2+|v|^2}}^\delta$ for $(u,v)\neq(0,0)$.
We have
\bel{Han-Lu-Zhu theorem}
\left\|\omega \S_\rho f\right\|_{\L^q(\R^{n+1})}~\leq~\C_{p~q}~\left\| f\sigma\right\|_{\L^p(\R^{2n+1})},\qquad 1<p\leq q<\infty
\eeq
if 
\bel{Formula}
 \gamma<{2n\over q},\qquad \delta<2n\left({p-1\over p}\right),\qquad \gamma+\delta\ge0,\qquad {\rho\over n+1}~=~{1\over p}-{1\over q}+{\gamma+\delta\over 2n+2}.
\eeq}
\begin{remark} Note that the two power weights $\omega,\sigma$  are defined in the subspace $\R^{2n}$. An analogue  two-weight $\L^p\mt\L^q$-norm inequality with   
\[\omega(u,v,t)~=~{\sqrt{|u|^2+|v|^2+|t|}}^{-\gamma}, \qquad \sigma(u,v,t)~=~{\sqrt{|u|^2+|v|^2+|t|}}^\delta\]
 can be found in the paper of  Han, Lu and Zhu \cite{Han-Lu-Zhu}.
\end{remark}

The proof of {\bf Han-Lu-Zhu theorem} was accomplished by using the language of fractional integrals defined in homogeneous spaces. In this paper, we first show that the  constraints inside (\ref{Formula}) are also necessary conditions for the $\L^p\mt\L^q$-norm inequality in (\ref{Han-Lu-Zhu theorem}). Conversely, we  give a new proof of  (\ref{Formula})  implying (\ref{Han-Lu-Zhu theorem}) for $1<p<q<\infty$ with a more direct approach.
\v
{\bf Theorem One}~~{\it Let $\S_\rho$ defined in (\ref{S_a})-(\ref{Omega^a}) for $0<\rho<n+1$. Suppose $\gamma,\delta\in\R$ and $\omega(u,v)={\sqrt{|u|^2+|v|^2}}^{-\gamma}, \sigma(u,v)={\sqrt{|u|^2+|v|^2}}^\delta$ for $(u,v)\neq(0,0)$.
We have
\bel{Result One}
\left\|\omega \S_\rho f\right\|_{\L^q(\R^{n+1})}~\leq~\C_{p~q}~\left\| f\sigma\right\|_{\L^p(\R^{2n+1})},\qquad 1<p< q<\infty
\eeq
if and only if
\bel{Formula One}
 \gamma<{2n\over q},\qquad \delta<2n\left({p-1\over p}\right),\qquad \gamma+\delta\ge0,\qquad {\rho\over n+1}~=~{1\over p}-{1\over q}+{\gamma+\delta\over 2n+2}.
\eeq}

Next, we extend {\bf Theorem One} to a multi-parameter setting  by replacing $\Omega^\rho$ with a larger kernel having singularity on every coordinate subspace. 

Observe that
\bel{kernel compara}
\Omega^\rho(u,v,t)~\leq~\left[{1\over |u||v|+|t|}\right]^{n+1-\rho},\qquad \hbox{\small{$(u,t)\neq(0,0)$~~ or~~ $(v,t)\neq(0,0)$}}.
\eeq
Furthermore, we find
\bel{Kernel computa}
\begin{array}{lr}\ds
\left[{1\over |u||v|+|t|}\right]^{n+1-\rho}
~\approx~ \left[{1\over |u|^2|v|^2+t^2}\right]^{{n+1\over 2}-{\rho\over2}}
\\\\ \ds~~~~~~~~~~~~~~~~~~~~~~~~~~~~
~=~|u|^{{\rho\over 2}-{n+1\over 2}}|v|^{{\rho\over 2}-{n+1\over 2}}|t|^{{\rho\over 2}-{n+1\over 2}} \left[{|u||v||t|\over |u|^2|v|^2+t^2}\right]^{{n+1\over 2}-{\rho\over2}}
\\\\ \ds~~~~~~~~~~~~~~~~~~~~~~~~~~~~
~=~|u|^{\big[{\rho\over 2}+{n-1\over 2}\big]-n}|v|^{\big[{\rho\over 2}+{n-1\over 2}\big]-n}|t|^{\big[{\rho\over 2}-{n-1\over 2}\big]-1} \left[{|u||v|\over |t|}+{|t|\over|u||v|}\right]^{-\big[{n+1\over 2}-{\rho\over2}\big]}.
\end{array}
\eeq
Above estimates  lead us to the following assertion. Let $\alpha, \beta\in\R$ and $\vartheta\ge0$. $\V^{\alpha\beta\vartheta}$ is a distribution in $\R^{2n+1}$ agree with 
\bel{V}\begin{array}{lr}\ds
\V^{\alpha\beta\vartheta}(u,v,t)~=~|u|^{\alpha-n}|v|^{\alpha-n}|t|^{\beta-1} \Bigg[ {|u||v|\over |t|}+{|t|\over |u||v|}\Bigg]^{-\vartheta},
\qquad
\hbox{\small{$u\neq0, v\neq0$, $t\neq0$}}.
\end{array}
\eeq
Define 
\bel{I}
\begin{array}{lr}\ds
\I_{\alpha\beta\vartheta} f(u,v,t)~=~\iiint_{\R^{2n+1}} f\left(\xi,\eta,\tau\right) \V^{\alpha\beta\vartheta}\Big[(u,v,t)\odot(\xi,\eta,\tau)^{-1}\Big]d\xi d\eta d\tau.
\end{array}
\eeq
This fractional integral operator is associated with Zygmund dilation, whereas
\[
\V^{\alpha\beta\vartheta}\Big[(ru,sv,rst)\odot(r\xi,s\eta,rs\tau)^{-1}\Big]~=~r^{\alpha+\beta-n-1}s^{\alpha+\beta-n-1} \V^{\alpha\beta\vartheta}\Big[(u,v,t)\odot(\xi,\eta,\tau)^{-1}\Big],\qquad r,s>0.\]

Singular integral operators with kernels having certain multi-parameter structures defined on Heisenberg group have been systematically studied, for example  by Phong and Stein \cite{Phong-Stein}, Ricci and Stein \cite{Ricci-Stein} and M\"{u}ller, Ricci and Stein \cite{Muller-Ricci-Stein}. Much less is known in this direction for fractional integration.
\v
{\bf Theorem Two}~~{\it Let $\I_{\alpha\beta\vartheta}$ defined in  (\ref{V})-(\ref{I}) for $\alpha,\beta\in\R$ and $\vartheta\ge0$. Suppose $\gamma,\delta\in\R$ and $\omega(u,v)={\sqrt{|u|^2+|v|^2}}^{-\gamma}, \sigma(u,v)={\sqrt{|u|^2+|v|^2}}^\delta$ for $(u,v)\neq(0,0)$.
We have
\bel{Result Two}
\left\| \omega\I_{\alpha\beta\vartheta} f\right\|_{\L^q(\R^{2n+1})}~\leq~\C_{p~q~\gamma~\delta}~\left\| f\sigma\right\|_{\L^p(\R^{2n+1})},\qquad 1<p< q<\infty
\eeq
if and only if
\bel{Formula Two}
\begin{array}{cc}\ds
 \gamma<{2n\over q},\qquad \delta<2n\left({p-1\over p}\right),\qquad\gamma+\delta\ge0,\qquad 
{\alpha+\beta\over n+1}~=~{1\over p}-{1\over q}+{\gamma+\delta\over 2n+2}
\\\\ \ds
 \vartheta~\ge~ \left| {\alpha-n\beta\over n+1}-{\gamma+\delta\over 2n+2}\right|;
\\\\ \ds
 n\left[{\alpha+\beta\over n+1}\right]+{\gamma+\delta\over 2n+2}-{n\over p}~<~\delta \qquad\hbox{for}\qquad \gamma\ge0,~\delta\leq0;
\\\\ \ds
 n\left[{\alpha+\beta\over n+1}\right]+{\gamma+\delta\over 2n+2}-n\left({q-1\over q}\right)~<~\gamma \qquad\hbox{for}\qquad \gamma\leq0,~\delta\ge0.
 \end{array}
\eeq}
\begin{remark} Recall (\ref{kernel compara})-(\ref{Kernel computa}).
By taking into account $\alpha={\rho\over 2}+{n-1\over 2}$, $\beta={\rho\over 2}-{n-1\over 2}$ and $ \vartheta^*= \left| {\alpha-n\beta\over n+1}-{\gamma+\delta\over 2n+2}\right|$
for $\rho, \gamma,\delta,p,q$ satisfying (\ref{Formula One}). We find
\[ \left[{1\over |u||v|+|t|}\right]^{n+1-\rho}~\lesssim~\V^{\alpha\beta\vartheta^*}(u,v,t),\qquad \hbox{\small{$u\neq0, v\neq0$, $t\neq0$}}.\]
This is equivalent to verify $\vartheta^*\leq {n+1\over 2}-{\rho\over 2}$. We omit the regarding computations.
\end{remark}
The remaining paper is organized as follows. First,  we prove {\bf Theorem One} in section 2. Section 3-5 are devoted to the proof of {\bf Theorem Two}.
In section 3, we show (\ref{Result Two}) implying (\ref{Formula Two}). In section 4, after a reformulation of $\I_{\alpha\beta\vartheta}$, we shall see that in the one-weight case, i.e : $\omega=\sigma$ occurred at $\gamma+\delta=0$, the $\L^p\mt\L^q$-norm inequality in (\ref{Result Two}) can be obtained by using an iteration argument. In contrast, this idea of iteration does not apply to  $\omega\neq\sigma$ whenever $\gamma+\delta>0$. In section 5, we develop a new framework to handle this two-weight case where the product space $\R^n\times\R^n$ is decomposed into an infinitely many dyadic cones. Each partial operator is defined on one of these dyadic cones. Essentially, it is a classical one-parameter fractional integral operator, satisfying the desired regularity. Moreover, its operator's norm decays as the eccentricity of the cone getting large.

\section{Proof of Theorem One}
\setcounter{equation}{0}

Let $\omega(u,v)={\sqrt{|u|^2+|v|^2}}^{-\gamma}$ and $\sigma(u,v)={\sqrt{|u|^2+|v|^2}}^\delta$ for $\gamma,\delta\in\R$ and $(u,v)\neq(0,0)$. 

Because $\S_\rho$ defined in (\ref{S_a})-(\ref{Omega^a}) for $0<\rho<n+1$  is self-adjoint, it is essential to have $\omega^q, \sigma^{-{p\over p-1}}$ locally integrable in $\R^{2n}$ for the $\L^p\mt\L^q$-norm inequality  in (\ref{Han-Lu-Zhu theorem}). Therefore, 
\bel{local inte A}
\gamma<{2n\over q},\qquad \delta<2n\left({p-1\over p}\right)
\eeq
are necessities.

Denote $\Q\subset\R^{2n}$ to be a cube parallel to the coordinates and $I\subset\R$ to be an interval. 
Consider
 \bel{f A}
 \begin{array}{lr}\ds
 f(u,v,t)~=~\sigma^{-{p\over p-1}}(u,v)\chi_{\Q\times I}(u,v,t)
 ~=~\sigma^{-{p\over p-1}}(u,v)\chi_{\Q}(u,v)\chi_{I}(t),\qquad \hbox{\small{$(u,v)\neq(0,0)$}}.
 \end{array}
 \eeq
Let $\vol\{\Q\}^{1\over n}=\vol\{I\}$. By changing variable $\tau\mt \tau-\mu\big(u\cdot\eta-v\cdot\xi\big)$ 
in $\left\|\omega\S_\rho f\right\|_{\L^q(\R^{2n+1})}$, we have
\bel{f=chi cone-para A}
\begin{array}{lr}\ds
\left\{\iiint_{\R^{2n+1}}\omega^q(u,v)\left\{\iiint_{\R^{2n+1}} f\left(\xi,\eta,\tau-\mu\big(u\cdot\eta-v\cdot\xi\big)\right) 
\Omega^\rho(u-\xi,v-\eta,t-\tau)
d\xi d\eta d\tau\right\}^q dudvdt\right\}^{1\over q}
\\\\ \ds
\geq~\Bigg\{\iiint_{\Q\times I}\omega^q(u,v)\Bigg\{\iiint_{\Q\times\R} \sigma^{-\frac{p}{p-1}}(\xi,\eta)\chi_{I}(\tau-\mu(u\cdot\eta-v\cdot\xi))
\\ \ds~~~~~~~~~~~~~~~~~~~~~~~~~~~~~~~~~~
\left[{1\over |u-\xi|^2+|v-\eta|^2+|t-\tau|}\right]^{n+1-\rho}
d\xi d\eta d\tau\Bigg\}^{q}du dv dt\Bigg\}^{\frac{1}{q}}
\\\\ \ds
\geq~\vol\{\Q\}^{{\rho\over n}-{n+1\over n}}
\left\{\iiint_{\Q\times I}\omega^{q}(u,v)\left\{\iint_\Q\sigma^{-\frac{p}{p-1}}(\xi,\eta) \left\{\int_{I-\mu(u\cdot\eta-v\cdot\xi)} d\tau\right\} d\xi d\eta \right\}^q dudvdt\right\}^{\frac{1}{q}}
\\\\ \ds
=~\vol\{\Q\}^{{\rho\over n}-{n+1\over n}+\big[1+{1\over q}\big]{1\over n}}
\left\{\iint_{\Q}\omega^{q}(u,v)dudv\right\}^{\frac{1}{q}}\iint_{\Q}\sigma^{-\frac{p}{p-1}}(\xi,\eta)d\xi d\eta.
\end{array}
\eeq
The $\L^p\mt\L^q$-norm inequality in (\ref{Han-Lu-Zhu theorem}) implies
\bel{chi norm ineq A}
\begin{array}{lr}\ds
\vol\{\Q\}^{{\rho\over n}-{n+1\over n}+\big[1+{1\over q}\big]{1\over n}}
\left\{\iint_{\Q}\omega^{q}(u,v)dudv\right\}^{\frac{1}{q}}\iint_{\Q}\sigma^{-\frac{p}{p-1}}(u,v)du dv
\\\\ \ds
~\leq~\C_{p~q}\left\{\iiint_{\Q\times I}\sigma^{-\frac{p}{p-1}}(u,v)dudvdt\right\}^{\frac{1}{p}}
~=~\C_{p~q}~\vol\{\Q\}^{{1\over n}{1\over p}} \left\{\iint_{\Q}\sigma^{-\frac{p}{p-1}}(u,v)dudv\right\}^{\frac{1}{p}}.
\end{array}
\eeq
From (\ref{f=chi cone-para A})-(\ref{chi norm ineq A}), we find
\bel{R R' A}
\begin{array}{lr}\ds
\vol\{\Q\}^{{\rho\over n}-{n+1\over n}+\big[1+{1\over q}-{1\over p}\big]{1\over n}}
\left\{\iint_{\Q}\omega^{q}(u,v)dudv\right\}^{\frac{1}{q}}\left\{\iint_{\Q}\sigma^{-\frac{p}{p-1}}(u,v)du dv\right\}^{p-1\over p}~=~
\\\\ \ds
\vol\{\Q\}^{\big[{\rho\over n+1}-{1\over p}+{1\over q}\big]{n+1\over n}}
\left\{{1\over \vol\{\Q\}}\iint_{\Q}\omega^{q}(u,v)dudv\right\}^{\frac{1}{q}}\left\{{1\over \vol\{\Q\}}\iint_{\Q}\sigma^{-\frac{p}{p-1}}(u,v)du dv\right\}^{p-1\over p}<\infty
\end{array}
\eeq
for every $\Q\subset\R^{2n}$.

A standard exercise of changing one-parameter dilation in (\ref{R R' A}) shows that
\bel{A homo}
{\rho\over n+1}~=~{1\over p}-{1\over q}+{\gamma+\delta\over 2n+2}
\eeq
 is an necessary homogeneity condition.

Let $\Q$ shrink to some $(u,v)\in\Q$ with $(u,v)\neq(0,0)$ inside (\ref{R R' A}). We have
\bel{limit Q}
\lim_{\vol\{\Q\}\mt0} \vol\{\Q\}^{\big[{\rho\over n+1}-{1\over p}+{1\over q}\big]{n+1\over n}} \omega(u,v)\sigma^{-1}(u,v)
\eeq
by applying Lebesgue differentiation theorem.
In order to have this limit finite, we need
\bel{Ineq > A}
\frac{\rho}{n+1}~\ge~\frac{1}{p}-\frac{1}{q}.
\eeq
By putting together  (\ref{A homo}) and (\ref{Ineq > A}), we find
\bel{gamma+delta}
\gamma+\delta~\ge~0.
\eeq
Recall $\Omega^\rho(u,v,t)$ defined in (\ref{Omega^a}). We have
\bel{Omega^a rewrite}
\begin{array}{lr}\ds
\Omega^\rho(u,v,t)~=~\left[{1\over |u|^2+|v|^2+|t|}\right]^{n+1-\rho}
\\\\ \ds~~~~~~~~~~~~~~~~
~=~\left[{1\over |u|^2+|v|^2+|t|}\right]^{n-\big({n\over n+1}\big)\rho-{\gamma+\delta\over 2n+2}+1-{\rho\over n+1}+{\gamma+\delta\over 2n+2}}
\\\\ \ds~~~~~~~~~~~~~~~~
~\leq~\left[{1\over |u|^2+|v|^2}\right]^{n-n\big[{\rho\over n+1}+{1\over n}{\gamma+\delta\over 2n+2}\big]}  |t|^{\big[{\rho\over n+1}-{\gamma+\delta\over 2n+2}\big]-1},
\qquad \hbox{\small{$(u,v)\neq(0,0)$,~ $t\neq0$}}.
\end{array}
\eeq
Note that a direct computation shows
\bel{homogeneity computa}
\begin{array}{lr}\ds
{\rho\over n+1}+{1\over n}{\gamma+\delta\over 2n+2}~=~{1\over p}-{1\over q}+{\gamma+\delta\over 2n}\qquad\hbox{\small{( ${\rho\over n+1}={1\over p}-{1\over q}+{\gamma+\delta\over 2n+2}$ )}}
\\\\ \ds~~~~~~~~~~~~~~~~~~~~~~~~~~~~
~<~1\qquad \hbox{\small{because~$\gamma<{2n\over q}$,~$\delta<2n\left({p-1\over p}\right)$}}.
\end{array}
\eeq
Let $\S_\rho$ defined in (\ref{S_a})-(\ref{Omega^a}) for $0<\rho<n+1$. 
By changing variable $\tau\mt \tau-\mu\big(u\cdot\eta-v\cdot\xi\big)$, we find
 \bel{S < =}
\begin{array}{lr}\ds
\S_\rho f(u,v,t)~=~\iiint_{\R^{2n+1}} f\left(\xi,\eta,\tau-\mu\big(u\cdot\eta-v\cdot\xi\big)\right) 
\Omega^\rho(u-\xi,v-\eta,t-\tau)
d\xi d\eta d\tau
\\\\ \ds~~~~~~~~~~~~~~~~~
~\leq~
\iiint_{\R^{2n+1}}  f\left(\xi,\eta,\tau-\mu\big(u\cdot\eta-v\cdot\xi\big)\right) 
\\ \ds~~~~~~~~~~~~~~~~~~~~~~~~~
\left[{1\over |u-\xi|^2+|v-\eta|^2}\right]^{n-n\big[{\rho\over n+1}+{1\over n}{\gamma+\delta\over 2n+2}\big]}  |t-\tau|^{\big[{\rho\over n+1}-{\gamma+\delta\over 2n+2}\big]-1}d\xi d\eta d\tau\qquad\hbox{\small{by (\ref{Omega^a rewrite})}}
\\\\ \ds~~~~~~~~~~~~~~~~~
~\doteq~\iint_{\R^{2n}} \left[{1\over |u-\xi|^2+|v-\eta|^2}\right]^{n-n\big[{\rho\over n+1}+{1\over n}{\gamma+\delta\over 2n+2}\big]} \F_{\rho\gamma\delta}(\xi,\eta, u,v,t)d\xi d\eta
\end{array}
\eeq
where
\bel{F'}
\F_{\rho\gamma\delta}(\xi,\eta, u,v,t)~=~\int_\R f\left(\xi,\eta,\tau-\mu\left(u \cdot\eta-v \cdot\xi\right)\right) |t-\tau|^{\big[{\rho\over n+1}-{\gamma+\delta\over 2n+2}\big]-1} d\tau.
\eeq
Because $\left[{1\over |u|^2+|v|^2}\right]^{n-n\big[{\rho\over n+1}+{1\over n}{\gamma+\delta\over 2n+2}\big]}$, $|t|^{\big[{\rho\over n+1}-{\gamma+\delta\over 2n+2}\big]-1}$ are positive definite, it is suffice to assert $f\ge0$.

Recall  {\bf Hardy-Littlewood-Sobolev theorem} and {\bf Stein-Weiss theorem} stated in the beginning of this paper. 
By applying (\ref{HLS Ineq}) with $\a={\rho\over n+1}-{\gamma+\delta\over 2n+2}={1\over p}-{1\over q}$ and $\N=1$, we obtain
\bel{F' regularity}
\begin{array}{lr}\ds
\left\{\int_{\R} \F^q_{\rho\gamma\delta}(\xi,\eta, u,v,t) dt \right\}^{1\over q}~\leq~\C_{p~q} \left\{\int_{\R} \Big[ f\left(\xi,\eta,t-\mu\left(u \cdot\eta-v \cdot\xi\right)\right)\Big]^p dt\right\}^{1\over p}
\\\\ \ds~~~~~~~~~~~~~~~~~~~~~~~~~~~~~~~~~~~~~~~~~~~
~=~\C_{p~q}~ \left\| f(\xi,\eta,\cdot)\right\|_{\L^p(\R)}
\end{array}
\eeq
regardless of $(u,v)\in\R^n\times\R^n$.

From (\ref{S < =})-(\ref{F'}), we have
\bel{Est One}
\begin{array}{lr}\ds
\left\{\iiint_{\R^{2n+1}}{\sqrt{|u|^2+|v|^2}}^{-\gamma q} \Big(\S_\rho f\Big)^q(u,v,t) dudvdt\right\}^{1\over q}
\\\\ \ds
~\leq~\Bigg\{ \iiint_{\R^{2n+1}}{\sqrt{|u|^2+|v|^2}}^{-\gamma q}
\\ \ds~~~~~~
\left\{\iint_{\R^{2n}} \left[{1\over |u-\xi|^2+|v-\eta|^2}\right]^{n-n\big[{\rho\over n+1}+{1\over n}{\gamma+\delta\over 2n+2}\big]} \F_{\rho\gamma\delta}(\xi,\eta, u,v,t)d\xi d\eta\right\}^q dudvdt\Bigg\}^{1\over q}
\\\\ \ds
\leq~      \Bigg\{\iint_{\R^{2n}}   {\sqrt{|u|^2+|v|^2}}^{-\gamma q}  
\\ \ds~~~~
 \left\{ \iint_{\R^{2n}}   \left[{1\over |u-\xi|^2+|v-\eta|^2}\right]^{n-n\big[{\rho\over n+1}+{1\over n}{\gamma+\delta\over 2n+2}\big]}  \left\{\int_\R \F^q_{\rho\gamma\delta}(\xi,\eta, u,v,t) dt \right\}^{1\over q} d\xi d\eta\right\}^q dudv\Bigg\}^{1\over q}
\\ \ds~~~~~~~~~~~~~~~~~~~~~~~~~~~~~~~~~~~~~~~~~~~~~~~~~~~~~~~~~~~~~~~~~~~~~~~~~~~~~~~~~~~~~
\hbox{\small{by Minkowski integral inequality}}
\\\\ \ds
\leq~\C_{p~q} \Bigg\{\iint_{\R^{2n}}  {\sqrt{|u|^2+|v|^2}}^{-\gamma q}  
\\ \ds~~~~~~~~~~~
  \left\{ \iint_{\R^{2n}}   \left[{1\over |u-\xi|^2+|v-\eta|^2}\right]^{n-n\big[{\rho\over n+1}+{1\over n}{\gamma+\delta\over 2n+2}\big]}    \left\| f(\xi,\eta,\cdot)\right\|_{\L^p(\R)} d\xi d\eta\right\}^q dudv\Bigg\}^{1\over q}~~~~~\hbox{\small{by (\ref{F' regularity})}}
  \\\\ \ds
 ~\leq~\C_{p~q~\gamma~\delta}   \left\{ \iint_{\R^{2n}}      \left\| f(\xi,\eta,\cdot)\right\|_{\L^p(\R)}^p\left[\sqrt{|\xi|^2+|\eta|^2}\right]^{p\delta} d\xi d\eta\right\}^{1\over p}  
\\ \ds~~~~~~~~~~~~~~~~~~~~~
 \hbox{\small{by (\ref{homogeneity computa}) and applying (\ref{Stein-Weiss Ineq})-(\ref{formula}) with $\a=n\big[{\rho\over n+1}+{1\over n}{\gamma+\delta\over 2n+2}\big]$ and $\N=n$}}
 \\\\ \ds
 ~=~\C_{p~q~\gamma~\delta}   \left\{ \iiint_{\R^{2n+1}}      \Big[f(\xi,\eta,\tau)\Big]^p\left[\sqrt{|\xi|^2+|\eta|^2}\right]^{p\delta} d\xi d\eta d\tau\right\}^{1\over p}. 
\end{array}
\eeq

\section{Some necessary constraints}
\setcounter{equation}{0}
Recall $\I_{\alpha\beta\vartheta}$ defined in (\ref{V})-(\ref{I}). By changing variable $\tau\mt \tau-\mu\big(u\cdot\eta-v\cdot\xi\big)$, we find
\bel{I rewrite}
\begin{array}{lr}\ds
\I_{\alpha\beta\vartheta} f(u,v,t)~=~\iiint_{\R^{2n+1}} f\left(\xi,\eta,\tau-\mu\big(u\cdot\eta-v\cdot\xi\big)\right) 
\\ \ds~~~~~~~~~~~~~~~~~~~~~~~~~~~~~
|u-\xi|^{\alpha-n}|v-\eta|^{\alpha-n}|t-\tau|^{\beta-1} \Bigg[ {|u-\xi||v-\eta|\over |t-\tau|}+{|t-\tau|\over |u-\xi||v-\eta|}\Bigg]^{-\vartheta}d\xi d\eta d\tau.
\end{array}
\eeq
Let $\omega(u,v)={\sqrt{|u|^2+|v|^2}}^{-\gamma}$ and $\sigma(u,v)={\sqrt{|u|^2+|v|^2}}^\delta$ for $\gamma,\delta\in\R$ and $(u,v)\neq(0,0)$.

By changing dilations $(u,v,t)\mt (r u, r v, r^2\lambda t)$ and $(\xi,\eta,\tau)\mt (r \xi, r \eta, r^2\lambda \tau)$ for $r>0$ and  $0<\lambda<1$ or $\lambda>1$,  we have
\bel{Dila I general}
\begin{array}{lr}\ds
\left\{ \iiint_{\R^{2n+1}} \omega^q(u,v)
\left\{\iiint_{\R^{2n+1}} f\left[r^{-1} \xi,r^{-1}\eta,r^{-2}\lambda^{-1}\big[\tau-\mu\lambda\big(u\cdot\eta-v\cdot\xi\big)\big]\right] \right.\right.
\\ \ds
\left. \left.|u-\xi|^{\alpha-n}|v-\eta|^{\alpha-n}|t-\tau|^{\beta-1} \Bigg[ {|u-\xi||v-\eta|\over |t-\tau|}+{|t-\tau|\over |u-\xi||v-\eta|}\Bigg]^{-\vartheta}
 d\xi d\eta d\tau\right\}^q dudvdt\right\}^{1\over q}
\\\\ \ds
=~r^{2\alpha+2\beta}r^{-\gamma} r^{2n+2\over q}\lambda^\beta\lambda^{1\over q}~\left\{ \iiint_{\R^{2n+1}} \Big[\sqrt{ |u|^2+|v|^2}\Big]^{-\gamma q}
\left\{\iiint_{\R^{2n+1}} f\left( \xi,\eta,\tau-\mu\big(u\cdot\eta-v\cdot\xi\big)\right) \right.\right.
\\ \ds~~~
\left. \left.|u-\xi|^{\alpha-n}|v-\eta|^{\alpha-n}|t-\tau|^{\beta-1} \Bigg[ {|u-\xi||v-\eta|\over \lambda|t-\tau|}+{\lambda|t-\tau|\over |u-\xi||v-\eta|}\Bigg]^{-\vartheta}
 d\xi d\eta d\tau\right\}^q dudvdt\right\}^{1\over q}
\\\\ \ds
\ge~r^{2\alpha+2\beta}r^{-\gamma}r^{2n+2\over q}\lambda^\beta\lambda^{1\over q}\left\{ \begin{array}{lr}\ds \lambda^\vartheta,\qquad0<\lambda<1,
\\ \ds
\lambda^{-\vartheta},\qquad \lambda>1
\end{array}\right.
\\ \ds~~~
\left\{ \iiint_{\R^{2n+1}}\Big[\sqrt{ |u|^2+|v|^2}\Big]^{-\gamma q}
\left\{\iiint_{\R^{2n+1}} f\left( \xi,\eta,\tau-\mu\big(u\cdot\eta-v\cdot\xi\big)\right) \right.\right.
\\ \ds~~~
\left. \left.|u-\xi|^{\alpha-n}|v-\eta|^{\alpha-n}|t-\tau|^{\beta-1} \Bigg[ {|u-\xi||v-\eta|\over |t-\tau|}+{|t-\tau|\over |u-\xi||v-\eta|}\Bigg]^{-\vartheta}
 d\xi d\eta d\tau\right\}^q dudvdt\right\}^{1\over q}.
\end{array}
\eeq
The $\L^p\mt\L^q$-norm inequality in (\ref{Result One}) implies that the last line of (\ref{Dila I general}) is bounded by
\bel{I L^p}
\begin{array}{lr}\ds
\left\{ \iiint_{\R^{2n+1}} \Big| f\left(r^{-1} \xi,r^{-1}\eta,r^{-2}\lambda^{-1} \tau\right) \Big|^p \left[\sqrt{|\xi|^2+|\eta|^2}\right]^{\delta p} d\xi d\eta d\tau\right\}^{1\over p}
\\\\ \ds
~=~r^{2n+2\over p} r^\delta\lambda^{1\over p} \left\| f\sigma\right\|_{\L^p(\R^{2n+1})},
\qquad
\hbox{\small{$(\xi,\eta,\tau)\mt(r\xi,r\eta, r^2\lambda \tau)$}}.
\end{array}
\eeq
This must be true for every $r>0$ and $0<\lambda<1$ or $\lambda>1$. We necessarily have
\bel{homogeneity}
{\alpha+\beta\over n+1}~=~{1\over p}-{1\over q}+{\gamma+\delta\over 2n+2}
\eeq
and
\bel{Constraints alphabeta'}
\beta+\vartheta~\ge~{1\over p}-{1\over q}\qquad\hbox{\small{or}}\qquad \beta-\vartheta~\leq~{1\over p}-{1\over q}.
\eeq
By adding  (\ref{homogeneity}) and (\ref{Constraints alphabeta'}) together, we find
\[
\begin{array}{lr}\ds
\vartheta~\ge~{n\beta-\alpha\over n+1}+{\gamma+\delta\over 2n+2}
\qquad\hbox{or}\qquad
\vartheta~\ge~{\alpha-n\beta\over n+1}-{\gamma+\delta\over 2n+2}.
\end{array}
\]
This further implies
\bel{theta>}
\vartheta~\ge~\left| {\alpha-n\beta\over n+1}-{\gamma+\delta\over 2n+2}\right|.
\eeq
Because $\I_{\alpha\beta\vartheta}$ is self-adjoint, it is essential to have $\omega^q, \sigma^{-{p\over p-1}}$ locally integrable. Therefore, 
\bel{local inte}
\gamma<{2n\over q},\qquad \delta<2n\left({p-1\over p}\right)
\eeq
are necessary.

Denote $\Rec=\Q_1\times\Q_2\times I\subset\R^n\times\R^n\times\R$  where $\Q_1,\Q_2$  are cubes in $\R^n$   parallel to the coordinates. Moreover, $I$ is an interval.  $\Rec'=\Q'_1\times\Q'_2\times I'$ is a translation of $\Rec$ defined as
\bel{R'}
\Rec'~=~\left\{(u,v,t)\colon \begin{array}{cc}\ds u_i=\xi_i+2\vol\{\Q_1\}^{1\over n},~v_i=\eta_i+2 \vol\{\Q_2\}^{1\over n},~i=1,2,\ldots,n
\\ \ds
t=\tau+2 \vol\{I\}
\end{array}~(\xi,\eta,\tau)\in\Rec
\right\}.
\eeq
Consider
 \bel{f}
 f(u,v,t)~=~\sigma^{-{p\over p-1}}(u,v)\chi_{\Q_1\times\Q_2}(u,v)\chi_{I}(t),\qquad\hbox{\small{$(u,v)\neq(0,0)$}}
 \eeq
 where $\chi$ is an indicator function.

 Let 
$\vol\{I\}=\vol\{\Q_1\}^{1\over n}\vol\{\Q_2\}^{1\over n}$.
We have
\bel{f=chi}
\begin{array}{lr}\ds
\left\|\omega\I_{\alpha\beta\vartheta}f\right\|_{\L^q(\R^{2n+1})}~\geq~
\\\\ \ds
\Bigg\{\iiint_{\Rec'}\omega^{q}(u,v)\Bigg\{\iiint_{\Q_1\times\Q_2\times\R} \sigma^{-\frac{p}{p-1}}(\xi,\eta)\chi_{I}(\tau-\mu(u\cdot\eta-v\cdot\xi))\left|u-\xi\right|^{\alpha-n}\left|v-\eta\right|^{\alpha-n}\left|t-\tau\right|^{\beta-1}
\\ \ds~~~~~~~~~~~~~~~~~~~~~~~~~~~~~~
\left[\frac{\left|u-\xi\right|\left|v-\eta\right|}{\left|t-\tau\right|}+\frac{\left|t-\tau\right|}{\left|u-\xi\right|\left|v-\eta\right|}\right]^{-\vartheta}d\xi d\eta d\tau\Bigg\}^{q}du dv dt\Bigg\}^{\frac{1}{q}}
\\\\ \ds
\geq~\vol\{\Q_1\}^{\frac{\alpha}{n}-1}\vol\{\Q_2\}^{\frac{\alpha}{n}-1}\vol\{I\}^{\beta-1}
\\ \ds~~
\left\{\iiint_{\Q_1'\times\Q_2'\times I'}\omega^{q}(u,v)\left\{\iint_{\Q_1\times\Q_2}\sigma^{-\frac{p}{p-1}}(\xi,\eta) \left\{\int_{I-\mu(u\cdot\eta-v\cdot\xi)} d\tau\right\} d\xi d\eta \right\}^q dudvdt\right\}^{\frac{1}{q}}
\\\\ \ds
=~\vol\{\Q_1\}^{\frac{\alpha}{n}-1}\vol\{\Q_2\}^{\frac{\alpha}{n}-1}\vol\{I\}^{\beta-1+\frac{1}{q}+1}
\\ \ds~~~
\left\{\iint_{\Q'_1\times\Q'_2}\omega^{q}(u,v)dudv\right\}^{\frac{1}{q}}\iint_{\Q_1\times\Q_2}\sigma^{-\frac{p}{p-1}}(\xi,\eta)d\xi d\eta
\\\\ \ds
=~\vol\{\Q_1\}^{\frac{\alpha}{n}-1}\vol\{\Q_2\}^{\frac{\alpha}{n}-1}\vol\{I\}^{\beta+\frac{1}{q}} \left\{\iint_{\Q'_1\times\Q'_2}\omega^{q}(u,v)dudv\right\}^{\frac{1}{q}}\iint_{\Q_1\times\Q_2}\sigma^{-\frac{p}{p-1}}(u,v)dudv.
\end{array}
\eeq
The norm inequality in (\ref{Result One}) implies
\bel{chi norm ineq}
\begin{array}{lr}\ds
\vol\{\Q_1\}^{\frac{\alpha}{n}-1}\vol\{\Q_2\}^{\frac{\alpha}{n}-1}\vol\{I\}^{\beta+\frac{1}{q}} \left\{\iint_{\Q'_1\times\Q'_2}\omega^{q}(u,v)dudv\right\}^{\frac{1}{q}}\iint_{\Q_1\times\Q_2}\sigma^{-\frac{p}{p-1}}(u,v)dudv
\\\\ \ds
~\leq~\C_{\alpha~\beta~p~q}~\vol\{I\}^{1\over p} \left\{\iint_{\Q_1\times\Q_2}\sigma^{-\frac{p}{p-1}}(u,v)dudv\right\}^{\frac{1}{p}}.
\end{array}
\eeq
By  taking into account $\vol\{I\}=\vol\{\Q_1\}^{1\over n}\vol\{\Q_2\}^{1\over n}$, we find
\bel{R R'}
\begin{array}{lr}\ds
~~~~~~~\vol\{\Q_1\}^{\frac{\alpha}{n}-1}\vol\{\Q_2\}^{\frac{\alpha}{n}-1}\vol\{I\}^{\beta+\frac{1}{q}-\frac{1}{p}}\left\{\iint_{\Q'_1\times\Q'_2}\omega^{q}(u,v)dudv\right\}^{\frac{1}{q}}\left\{\iint_{\Q_1\times\Q_2}\sigma^{-\frac{p}{p-1}}(u,v)dudv\right\}^{\frac{p-1}{p}}
\\\\ \ds
~=~\vol\{\Q_1\}^{\big[\frac{\alpha+\beta}{n+1}-\big(\frac{1}{p}-\frac{1}{q}\big)\big]\frac{n+1}{n}}
\vol\{\Q_2\}^{\big[\frac{\alpha+\beta}{n+1}-\big(\frac{1}{p}-\frac{1}{q}\big)\big]\frac{n+1}{n}}
\\ \ds~~~~~
\left\{\frac{1}{\vol\{\Q'_1\}\vol\{\Q'_2\}}\iint_{\Q'_1\times\Q'_2}\omega^q(u,v)dudv\right\}^{\frac{1}{q}}\left\{\frac{1}{\vol\{\Q_1\}\vol\{\Q_2\}}\iint_{\Q_1\times\Q_2}\sigma^{-\frac{p}{p-1}}(u,v)dudv\right\}^{\frac{p-1}{p}}
\\ \ds~~~~~~~~~~~~~~~~~~~~~~~~~~~~~~~~~~~~~~~~~~~~~~~~~~~~~~~~~~~~~~~~~~~~~~~~~~~~~~~~~~~~~~~~~~~~~~~~~~~~~~~~~~~~~~~~~~~~~~~~~~~~~~~~~~~~~~~~~~~~
~<~\infty
\end{array}
\eeq
for every $\Q_1\times\Q_2\subset\R^n\times\R^n$.

Note that (\ref{R R'}) holds for every $\Q_1\times\Q_2\subset\R^n\times\R^n$. Suppose $\Q_2$ centered on the origin and $\vol\{\Q_2\}^{1\over n}=1$.
Let $\Q_1$ shrink to $u\in\Q_1$. Simultaneously, as defined in (\ref{R'}), $\Q_1'$ shrinks to some $u'\in\Q_1'$ and $\vol\{\Q'_2\}^{1\over n}=1$.
By applying Lebesgue differentiation theorem, we find 
\bel{last line}
\begin{array}{lr}\ds
\lim_{\vol\{\Q_1\}\mt0}  \vol\{\Q_1\}^{\big[\frac{\alpha+\beta}{n+1}-\big(\frac{1}{p}-\frac{1}{q}\big)\big]\frac{n+1}{n}}
\left\{\int_{\Q'_2}\omega^q(u',v)dv\right\}^{\frac{1}{q}}\left\{\int_{\Q_2}\sigma^{-\frac{p}{p-1}}(u,v)dv\right\}^{\frac{p-1}{p}}
~<~\infty.
\end{array}
\eeq
Clearly, the product of two integral terms in (\ref{last line}) never vanishes. We must have
$\frac{\alpha+\beta}{n+1}\ge\frac{1}{p}-\frac{1}{q}$ in order to bound the limit as $\vol\{\Q_1\}\mt0$.
This together with the homogeneity condition in (\ref{homogeneity}) imply
\bel{gamma+delta}
\gamma+\delta~\ge~0.
\eeq
For brevity of computation, denote 
\bel{zeta}
\zeta~=~n\left[{\alpha+\beta\over n+1}\right]+{\gamma+\delta\over 2n+2}.
\eeq
We find
\bel{zeta homo}
\begin{array}{lr}\ds
\zeta~=~{n\over p}-{n\over q}+{\gamma+\delta\over 2}\qquad\hbox{\small{\big(~${\alpha+\beta\over n+1}={1\over p}-{1\over q}+{\gamma+\delta\over 2n+2}$~\big)}};
\\\\ \ds
0~<~\zeta~=~{n\over p}-{n\over q}+{\gamma+\delta\over 2}\qquad\hbox{\small{(~$\gamma+\delta\ge0$,~$1<p< q<\infty$)~}}
\\\\ \ds~
~<~{n\over p}-{n\over q}+{n\over q}+n\left({p-1\over p}\right)~=~n.\qquad\hbox{\small{$\big(~\gamma<{2n\over q}$,~$\delta<2n\big({p-1\over p}\big)~\big)$}}
\end{array}
\eeq
Moreover,  a direct computation shows
\bel{zeta equal}
\begin{array}{lr}\ds
\left[\frac{\alpha+\beta}{n+1}-\left(\frac{1}{p}-\frac{1}{q}\right)\right]\frac{n+1}{n}~=~\frac{\alpha+\beta}{n+1}-\left(\frac{1}{p}-\frac{1}{q}\right)+{1\over n}{\gamma+\delta\over 2n+2}\qquad\hbox{\small{by (\ref{homogeneity})}}
\\\\ \ds~~~~~~~~~~~~~~~~~~~~~~~~~~~~~~~~~~~~~~~~
~=~{\zeta\over n}-\left({1\over p}-{1\over q}\right).
\end{array}
\eeq
From (\ref{R R'}) and (\ref{zeta equal}), we obtain
\bel{Characteristic}
\begin{array}{lr}\ds
\sup_{\Q_1\times\Q_2\subset\R^n\times\R^n}~\vol\{\Q_1\}^{{\zeta\over n}-\big(\frac{1}{p}-\frac{1}{q}\big)}
\vol\{\Q_2\}^{{\zeta\over n}-\big(\frac{1}{p}-\frac{1}{q}\big)}
\\ \ds~~~~~~~~~~~~~~~~~~~~~
\left\{\frac{1}{\vol\{\Q'_1\}\vol\{\Q'_2\}}\iint_{\Q'_1\times\Q'_2}\left[\sqrt{|u|^2+|v|^2}\right]^{-\gamma q}dudv\right\}^{\frac{1}{q}}
\\ \ds~~~~~~~~~~~~~~~~~~~~~
\left\{\frac{1}{\vol\{\Q_1\}\vol\{\Q_2\}}\iint_{\Q_1\times\Q_2}\left[\sqrt{|u|^2+|v|^2}\right]^{-\delta\frac{p}{p-1}}dudv\right\}^{\frac{p-1}{p}}
~<~\infty.
\end{array}
\eeq

\subsection{ Case One: $\gamma\ge0$,  $\delta\leq0$}
Suppose $\gamma+\delta=0$. Let $\zeta$ defined in (\ref{zeta}). From (\ref{homogeneity}) and (\ref{zeta equal}), we find 
\bel{balance}
{\zeta\over n}~=~{1\over p}-{1\over q}.
\eeq 
Recall $\Rec'=\Q'_1\times\Q'_2\times I$ defined in (\ref{R'}) which is a translation of $\Rec=\Q_1\times\Q_2\times I$. 
We consider $\Q'_1\times\Q'_2$ centered on the origin of $\R^n\times\R^n$.
Let $\vol\{\Q_2\}^{1\over n}=\vol\{\Q'_2\}^{1\over n}=1$ and $\Q_1$ shrink to some $u\in\Q_1$ whereas $\Q'_1$ shrink to $0$. 

From (\ref{Characteristic})-(\ref{balance}), by applying Lebesgue differentiation theorem, we have
\bel{Characteristic1}
\begin{array}{lr}\ds
\left\{\frac{1}{\vol\{\Q'_2\}}\int_{\Q'_2}|v|^{-\gamma q}dv\right\}^{\frac{1}{q}}
\left\{\frac{1}{\vol\{\Q_2\}}\int_{\Q_2}\Big[\sqrt{|u|^2+|v|^2}\Big]^{-\delta\frac{p}{p-1}}dv\right\}^{\frac{p-1}{p}}
~<~\infty
\end{array}
\eeq
for every $\Q_2\subset\R^n$. This suggests 
\bel{gamma<nq}
\gamma~<~{n\over q}\qquad\Longrightarrow\qquad \zeta-{n\over p}~=~-{n\over q}<-\gamma~=~\delta
\eeq
as an necessity.

Suppose $\gamma+\delta>0$. From (\ref{homogeneity}) and (\ref{zeta equal}), we find 
\bel{subbalance}
{\zeta\over n}~>~{1\over p}-{1\over q}.
\eeq 
For every $\Q_1\times\Q_2\subset\R^n\times\R^n$, we define
\bel{A Q}
\begin{array}{lr}\ds
\A_{p~q}^{\zeta~\gamma~\delta}(\Q_1\times\Q_2)~=~\vol\{\Q_1\}^{{\zeta\over n}-\big(\frac{1}{p}-\frac{1}{q}\big)}
\vol\{\Q_2\}^{{\zeta\over n}-\big(\frac{1}{p}-\frac{1}{q}\big)} 
\\ \ds~~~~~~~~~~~~~~~~~~~~~~~~~~~~~~~~~
\left\{\frac{1}{\vol\{\Q'_1\}\vol\{\Q'_2\}}\iint_{\Q'_1\times\Q'_2}\left[\sqrt{|u|^2+|v|^2}\right]^{-\gamma q}dudv\right\}^{\frac{1}{q}}
\\ \ds~~~~~~~~~~~~~~~~~~~~~~~~~~~~~~~~~
\left\{\frac{1}{\vol\{\Q_1\}\vol\{\Q_2\}}\iint_{\Q_1\times\Q_2}\left[\sqrt{|u|^2+|v|^2}\right]^{-\delta\frac{p}{p-1}}dudv\right\}^{\frac{p-1}{p}}.
\end{array}
\eeq
Moreover, denote
\bel{Q^k'}
{\Q'}_1^k~=~\Q'_1\cap\Big\{ 2^{-k-1}\leq|u|<2^{-k}\Big\},\qquad {\Q'}_2^k~=~\Q'_2\cap\Big\{ 2^{-k-1}\leq|v|<2^{-k}\Big\},\qquad k\ge0.
\eeq 
Let $\vol\{\Q_2\}^{1\over n}=\vol\{\Q'_2\}^{1\over n}=1$ and $\vol\{\Q_1\}^{1\over n}=\vol\{\Q'_1\}^{1\over n}=\lambda$ for $0<\lambda<1$. From (\ref{A Q})-(\ref{Q^k'}), we have
\bel{Chara rewrite}
\begin{array}{lr}\ds
\left[\A_{p~q}^{\zeta~\gamma~\delta}(\Q_1\times\Q_2)\right]^q~=~
\lambda^{q\big[\zeta-\big(\frac{n}{p}-\frac{n}{q}\big)\big]}
\left\{\frac{1}{\lambda^n}\iint_{\Q'_1\times\Q'_2}\left[\sqrt{|u|^2+|v|^2}\right]^{-\gamma q}dudv\right\}
\\ \ds~~~~~~~~~~~~~~~~~~~~~~~~~~~~~~~~~~~~~~~~~~~~~~~~~~~~~~~~~
\left\{\frac{1}{\lambda^n}\iint_{\Q_1\times\Q_2}\left[\sqrt{|u|^2+|v|^2}\right]^{-\delta\frac{p}{p-1}}dudv\right\}^{\big[\frac{p-1}{p}\big]q}
\\\\ \ds
~=~\lambda^{q\big[\zeta-\big(\frac{n}{p}-\frac{n}{q}\big)\big]}
\sum_{k\ge0}\left\{\frac{1}{\lambda^n}\iint_{\Q'_1\times{\Q'}_2^k}\left[\sqrt{|u|^2+|v|^2}\right]^{-\gamma q}dudv\right\}
\\ \ds~~~~~~~~~~~~~~~~~~~~~~~~~~~~~~~~
\left\{\frac{1}{\lambda^n}\iint_{\Q_1\times\Q_2}\left[\sqrt{|u|^2+|v|^2}\right]^{-\delta\frac{p}{p-1}}dudv\right\}^{\big[\frac{p-1}{p}\big]q}
\\\\ \ds
~\doteq~\sum_{k\ge0} \A_k(\lambda).
\end{array}
\eeq
Lebesgue's Differentiation Theorem implies
\bel{Chara EST}
\begin{array}{lr}\ds
\lim_{\lambda\mt0}~ \frac{1}{\lambda^n}\iint_{\Q'_1\times{\Q'}_2^k}\left[\sqrt{|u|^2+|v|^2}\right]^{-\gamma q}dudv
~=~ \int_{{\Q'}_2^k} |v|^{-\gamma q} dv.
\end{array}
\eeq
Because $\delta\leq0$ and $\zeta>{n\over p}-{n\over q}$, we find
\bel{A0}
\A_k(0)~=~0,\qquad k\ge0.
\eeq
Note that (\ref{A0}) is true if $\zeta-\big({n\over p}-{n\over q}\big)$ in (\ref{Chara rewrite}) is replaced by any  smaller positive number. 
Therefore,  $\A_k(\lambda)$ is H\"{o}lder continuous $w.r.t~\lambda$ whose exponent remains strictly positive as $k\mt\infty$. 
Recall (\ref{Characteristic}). We have $\sum_{k\ge0}\A_k(\lambda)\leq\C_{\alpha~\gamma~\delta~q}$ for every $\lambda>0$.  Consequently,   $\sum_{k\ge0}\A_k(\lambda)$ is continuous at $\lambda=0$ and
\bel{Chara Lim}
\lim_{\lambda\mt0} ~\sum_{k\ge0}\A_k(\lambda)~=~0.
\eeq
A direct computation shows
\bel{Characteristic Est Case1}
\begin{array}{lr}\ds
\left[\A_{p~q}^{\zeta~\gamma~\delta}(\Q_1\times\Q_2)\right]^q~=~
\lambda^{q\big[\zeta-\big(\frac{n}{p}-\frac{n}{q}\big)\big]}
\left\{\frac{1}{\lambda^n}\iint_{\Q'_1\times\Q'_2}\left[\sqrt{|u|^2+|v|^2}\right]^{-\gamma q}dudv\right\}
\\ \ds~~~~~~~~~~~~~~~~~~~~~~~~~~~~~~~~~~~~~~~~~~~~~~~~~~~~~~~~~
\left\{\frac{1}{\lambda^n}\iint_{\Q_1\times\Q_2}\left[\sqrt{|u|^2+|v|^2}\right]^{-\delta\frac{p}{p-1}}dudv\right\}^{\big[\frac{p-1}{p}\big]q}
\\\\ \ds
~\ge~\C\lambda^{q\big[\zeta-\big(\frac{n}{p}-\frac{n}{q}\big)\big]}
\int_{\Q'_2}\left[\sqrt{\lambda^2+|v|^2}\right]^{-\gamma q}dv\qquad\hbox{\small{ ( $\delta\leq0$,~~$\vol\{\Q_1\}^{1\over n}=\vol\{\Q'_1\}^{1\over n}=\lambda$ )}}
\\\\ \ds
~\ge~\C \lambda^{q\big[\zeta-\big(\frac{n}{p}-\frac{n}{q}\big)\big]}\int_{0<|v|\leq\lambda}\left({1\over \lambda}\right)^{\gamma q} dv
~=~\C_{\gamma~q}~\lambda^{n-\gamma q+q\big[\zeta-\big(\frac{n}{p}-\frac{n}{q}\big)\big]}.
\end{array}
\eeq
From (\ref{Chara Lim})-(\ref{Characteristic Est Case1}), by using $\zeta={n\over p}-{n\over q}+{\gamma+\delta\over2}$ as shown in (\ref{zeta homo}), we find
\bel{Constraint est.1}
\begin{array}{lr}\ds
{n\over q}-\gamma+\zeta-\left(\frac{n}{p}-\frac{n}{q}\right)~>~0\qquad\Longrightarrow
\\\\ \ds
\zeta~<~{n\over q}-\gamma+2\zeta-\left(\frac{n}{p}-\frac{n}{q}\right)~=~{n\over q}-\gamma+ \left(\frac{n}{p}-\frac{n}{q}\right)+\gamma+\delta
\\\\ \ds~~~~~~~~~~~~~~~~~~~~~~~~~~~~~~~~~~~~~~~~~~~~~~
~=~{n\over p}+\delta.
\end{array}
\eeq
Recall $\zeta=n\big[{\alpha+\beta\over n+1}\big]+{\gamma+\delta\over 2n+2}$. By putting together (\ref{gamma<nq}) and (\ref{Constraint est.1}), we obtain
\bel{Case One Constraint}
 n\left[{\alpha+\beta\over n+1}\right]+{\gamma+\delta\over 2n+2}-{n\over p}~<~\delta \qquad\hbox{for}\qquad \gamma\ge0,~\delta\leq0.
 \eeq

\subsection{ Case Two: $\gamma\leq0$,  $\delta\ge0$}
Suppose $\gamma+\delta=0$. From (\ref{homogeneity}) and (\ref{zeta equal}), we find 
${\zeta\over n}={1\over p}-{1\over q}={q-1\over q}-{p-1\over p}$ as shown in (\ref{balance}).
The estimate in (\ref{Characteristic1}) suggests
\bel{delta<np'}
\delta~<~n\left({p-1\over p}\right)\qquad\Longrightarrow\qquad \zeta-n\left({q-1\over q}\right)~=~-n\left({p-1\over p}\right)<-\delta~=~\gamma
\eeq
as an necessity.

Suppose $\gamma+\delta>0$. From (\ref{homogeneity}) and (\ref{zeta equal}), we find 
${\zeta\over n}>{1\over p}-{1\over q}$ as (\ref{subbalance}).

For every $\Q_1\times\Q_2\subset\R^n\times\R^n$, $\A_{p~q}^{\zeta~\gamma~\delta}(\Q_1\times\Q_2)$ is defined in (\ref{A Q}). Denote
\bel{Q^k}
\Q_1^k~=~\Q_1\cap\Big\{ 2^{-k-1}\leq|u|<2^{-k}\Big\},\qquad \Q_2^k~=~\Q_2\cap\Big\{ 2^{-k-1}\leq|v|<2^{-k}\Big\},\qquad k\ge0.
\eeq 
As before, suppose $\vol\{\Q_2\}^{1\over n}=\vol\{\Q'_2\}^{1\over n}=1$ and $\vol\{\Q_1\}^{1\over n}=\vol\{\Q'_1\}^{1\over n}=\lambda$ for $0<\lambda<1$. From (\ref{A Q}) and (\ref{Q^k}), we have
\bel{Chara rewrite again}
\begin{array}{lr}\ds
\left[\A_{p~q}^{\zeta~\gamma~\delta}(\Q_1\times\Q_2)\right]^{p\over p-1}~=~
\lambda^{{p\over p-1}\big[\zeta-\big(\frac{n}{p}-\frac{n}{q}\big)\big]}
\left\{\frac{1}{\lambda^n}\iint_{\Q'_1\times\Q'_2}\left[\sqrt{|u|^2+|v|^2}\right]^{-\gamma q}dudv\right\}^{{1\over q}{p\over p-1}}
\\ \ds~~~~~~~~~~~~~~~~~~~~~~~~~~~~~~~~~~~~~~~~~~~~~~~~~~~~~~~~~~~~~~~~
\left\{\frac{1}{\lambda^n}\iint_{\Q_1\times\Q_2}\left[\sqrt{|u|^2+|v|^2}\right]^{-\delta\frac{p}{p-1}}dudv\right\}
\\\\ \ds
~=~\lambda^{{p\over p-1}\big[\zeta-\big(\frac{n}{p}-\frac{n}{q}\big)\big]}
\left\{\frac{1}{\lambda^n}\iint_{\Q'_1\times\Q'_2}\left[\sqrt{|u|^2+|v|^2}\right]^{-\gamma q}dudv\right\}^{{1\over q}{p\over p-1}}
\\ \ds~~~~~~~~~~~~~~~~~~~~~~~
\sum_{k\ge0} \left\{\frac{1}{\lambda^n}\iint_{\Q_1\times{\Q}_2^k}\left[\sqrt{|u|^2+|v|^2}\right]^{-\delta\frac{p}{p-1}}dudv\right\}
\\\\ \ds
~\doteq~\sum_{k\ge0} \B_k(\lambda).
\end{array}
\eeq
Lebesgue's Differentiation Theorem implies
\bel{Chara EST again}
\begin{array}{lr}\ds
\lim_{\lambda\mt0}~ \frac{1}{\lambda^n}\iint_{\Q_1\times{\Q}_2^k}\left[\sqrt{|u|^2+|v|^2}\right]^{-\delta{p\over p-1}}dudv
~=~ \int_{{\Q}_2^k} |v|^{-\delta {p\over p-1}} dv.
\end{array}
\eeq
Because $\gamma\leq0$ and $\zeta>{n\over p}-{n\over q}$, we find
\bel{B0}
\B_k(0)~=~0,\qquad k\ge0.
\eeq
As same as (\ref{A0}),
the estimate in (\ref{B0})  is true if $\zeta-\big({n\over p}-{n\over q}\big)$ in (\ref{Chara rewrite again}) is replaced by a  smaller positive number. 
Therefore,  $\B_k(\lambda)$ is H\"{o}lder continuous $w.r.t~\lambda$ whose exponent remains strictly positive as $k\mt\infty$. 

Recall (\ref{Characteristic}). We have $\sum_{k\ge0}\B_k(\lambda)\leq\C_{\alpha~\gamma~\delta~q}$ for every $\lambda>0$.  Consequently,   $\sum_{k\ge0}\B_k(\lambda)$ is continuous at $\lambda=0$ and
\bel{Chara Lim again}
\lim_{\lambda\mt0} ~\sum_{k\ge0}\B_k(\lambda)~=~0.
\eeq
A direct computation shows
\bel{Characteristic Est Case2}
\begin{array}{lr}\ds
\left[\A_{p~q}^{\zeta~\gamma~\delta}(\Q_1\times\Q_2)\right]^{p\over p-1}~=~
\lambda^{{p\over p-1}\big[\zeta-\big(\frac{n}{p}-\frac{n}{q}\big)\big]}
\left\{\frac{1}{\lambda^n}\iint_{\Q'_1\times\Q'_2}\left[\sqrt{|u|^2+|v|^2}\right]^{-\gamma q}dudv\right\}^{{1\over q}{p\over p-1}}
\\ \ds~~~~~~~~~~~~~~~~~~~~~~~~~~~~~~~~~~~~~~~~~~~~~~~~~~~~~~~~~~~~~~~~
\left\{\frac{1}{\lambda^n}\iint_{\Q_1\times\Q_2}\left[\sqrt{|u|^2+|v|^2}\right]^{-\delta\frac{p}{p-1}}dudv\right\}
\\\\ \ds
~\ge~\C\lambda^{{p\over p-1}\big[\zeta-\big(\frac{n}{p}-\frac{n}{q}\big)\big]}
\int_{\Q_2}\left[\sqrt{\lambda^2+|v|^2}\right]^{-\delta {p\over p-1}}dv\qquad\hbox{\small{ ( $\gamma\leq0$,~~$\vol\{\Q_1\}^{1\over n}=\vol\{\Q'_1\}^{1\over n}=\lambda$ )}}
\\\\ \ds
~\ge~\C \lambda^{{p\over p-1}\big[\zeta-\big(\frac{n}{p}-\frac{n}{q}\big)\big]}\int_{0<|v|\leq\lambda}\left({1\over \lambda}\right)^{\delta {p\over p-1}} dv
~=~\C_{\delta~p}~\lambda^{n-\delta\big({p\over p-1}\big)+{p\over p-1}\big[\zeta-\big(\frac{n}{p}-\frac{n}{q}\big)\big]}.
\end{array}
\eeq
From (\ref{Chara Lim again})-(\ref{Characteristic Est Case2}), by using $\zeta={n\over p}-{n\over q}+{\gamma+\delta\over2}=n\big[{q-1\over q}-{p-1\over p}\big]+{\gamma+\delta\over 2}$ in (\ref{zeta homo}), we find
\bel{Constraint est.2}
\begin{array}{lr}\ds
n\left({p-1\over p}\right)-\delta+\zeta-\left(\frac{n}{p}-\frac{n}{q}\right)~>~0\qquad\Longrightarrow
\\\\ \ds
\zeta~<~n\left({p-1\over p}\right)-\delta+2\zeta-n\left[{q-1\over q}-{p-1\over p}\right]~=~n\left({p-1\over p}\right)-\delta+ n\left[{q-1\over q}-{p-1\over p}\right]+\gamma+\delta
\\\\ \ds~~~~~~~~~~~~~~~~~~~~~~~~~~~~~~~~~~~~~~~~~~~~~~~~~~~~~~~~~~~~~~~~~~~~~~~~~
~=~n\left({q-1\over q}\right)+\gamma.
\end{array}
\eeq
Recall $\zeta=n\big[{\alpha+\beta\over n+1}\big]+{\gamma+\delta\over 2n+2}$. By putting together (\ref{delta<np'}) and (\ref{Constraint est.2}), we obtain
\bel{Case Two Constraint}
 n\left[{\alpha+\beta\over n+1}\right]+{\gamma+\delta\over 2n+2}-n\left({q-1\over q}\right)~<~\gamma \qquad\hbox{for}\qquad \gamma\leq0,~\delta\ge0.
 \eeq

\section{Reformulation of $\I_{\alpha\beta\vartheta}$}
\setcounter{equation}{0}
Recall $\V^{\alpha\beta\vartheta}(u,v,t)$  defined in (\ref{V}) for $u\neq0,v\neq0,t\neq0$  and $\vartheta\ge\left| {\alpha-n\beta\over n+1}-{\gamma+\delta\over 2n+2}\right|$. 
Suppose  $2\alpha-2n\beta-\gamma-\delta\ge0 $. We have 
\bel{EST1}
\begin{array}{lr}\ds
\V^{\alpha\beta\vartheta}(u,v,t)~\leq~|u|^{\alpha-n}|v|^{\alpha-n}|t|^{\beta-1} \Bigg[ {|u||v|\over |t|}+{|t|\over |u||v|}\Bigg]^{-\big[{\alpha-n\beta\over n+1}-{\gamma+\delta\over 2n+2}\big]}
\\\\ \ds~~~~~~~~~~~~~~~~~~~~
~\leq~|u|^{\alpha-n}|v|^{\alpha-n}|t|^{\beta-1} \Bigg[ {|u||v|\over |t|}\Bigg]^{-\big[{\alpha-n\beta\over n+1}-{\gamma+\delta\over 2n+2}\big]}
\\\\ \ds~~~~~~~~~~~~~~~~~~~~
~=~|u|^{n\Big[{\alpha+\beta\over n+1}\Big]+{\gamma+\delta\over 2n+2}-n}|v|^{n\Big[{\alpha+\beta\over n+1}\Big]+{\gamma+\delta\over 2n+2}-n} |t|^{{\alpha+\beta\over n+1}-{\gamma+\delta\over 2n+2}-1},
\qquad\hbox{\small{ $u\neq0$, $v\neq 0$,   $t\neq0$}}.
\end{array}
\eeq
Suppose $2\alpha-2n\beta-\gamma-\delta\leq0 $. We find  
\bel{EST2}
\begin{array}{lr}\ds
\V^{\alpha\beta\vartheta}(u,v,t)~\leq~|u|^{\alpha-n}|v|^{\alpha-n}|t|^{\beta-1} \Bigg[ {|u||v|\over |t|}+{|t|\over |u||v|}\Bigg]^{{\alpha-n\beta\over n+1}-{\gamma+\delta\over 2n+2}}
\\\\ \ds~~~~~~~~~~~~~~~~~~~~
~\leq~|u|^{\alpha-n}|v|^{\alpha-n}|t|^{\beta-1} \Bigg[ {|t|\over |u||v|}\Bigg]^{{\alpha-n\beta\over n+1}-{\gamma+\delta\over 2n+2}}
\\\\ \ds~~~~~~~~~~~~~~~~~~~~
~=~|u|^{n\Big[{\alpha+\beta\over n+1}\Big]+{\gamma+\delta\over 2n+2}-n}|v|^{n\Big[{\alpha+\beta\over n+1}\Big]+{\gamma+\delta\over 2n+2}-n} |t|^{{\alpha+\beta\over n+1}-{\gamma+\delta\over 2n+2}-1}
,\qquad\hbox{\small{ $u\neq0$, $v\neq 0$,   $t\neq0$}}.
\end{array}
\eeq
As (\ref{zeta}), we write $\zeta=n\left[{\alpha+\beta\over n+1}\right]+{\gamma+\delta\over 2n+2}$ where $ 0<\zeta<n$. 
Let $\I_{\alpha\beta\vartheta}$ defined in (\ref{V})-(\ref{I}). From now on,  we assert $f\ge0$.
By changing variable $\tau\mt \tau-\mu\big(u\cdot\eta-v\cdot\xi\big)$, we have 
 \bel{I < =}
\begin{array}{lr}\ds
\I_{\alpha\beta\vartheta} f(u,v,t)~=~\iiint_{\R^{2n+1}} f\left(\xi,\eta,\tau-\mu\big(u\cdot\eta-v\cdot\xi\big)\right) \V^{\alpha\beta\vartheta}(u-\xi,v-\eta,t-\tau)d\xi d\eta d\tau
\\\\ \ds~~~~~~~~~~~~~~~~~~~~
~\leq~\iiint_{\R^{2n+1}} f\left(\xi,\eta,\tau-\mu\big(u\cdot\eta-v\cdot\xi\big)\right) 
\\ \ds~~~~~~~~~~~~~~~~~~~~~~~~~~~~~
|u-\xi|^{\alpha-n}|v-\eta|^{\alpha-n}|t-\tau|^{\beta-1} \Bigg[ {|u-\xi||v-\eta|\over |t-\tau|}+{|t-\tau|\over |u-\xi||v-\eta|}\Bigg]^{-\big|{\alpha-n\beta\over n+1}-{\gamma+\delta\over 2n+2}\big|}d\xi d\eta d\tau
\\\\ \ds~~~~~~~~~~~~~~~~~~~~
~\leq~
\iiint_{\R^{2n+1}}  f\left(\xi,\eta,\tau-\mu\big(u\cdot\eta-v\cdot\xi\big)\right) 
\\ \ds~~~~~~~~~~~~~~~~~~~~~~~~~~~~
|u-\xi|^{\zeta-n}|v-\eta|^{\zeta-n} |t-\tau|^{{\alpha+\beta\over n+1}-{\gamma+\delta\over 2n+2}-1}d\xi d\eta d\tau\qquad\hbox{\small{by (\ref{EST1})-(\ref{EST2})}}
\\\\ \ds~~~~~~~~~~~~~~~~~~~~
~\doteq~\iint_{\R^{2n}} |u-\xi|^{\zeta-n}|v-\eta|^{\zeta-n} \F_{\alpha\beta\gamma\delta}(\xi,\eta, u,v,t)d\xi d\eta
\end{array}
\eeq
where
\bel{F}
\F_{\alpha\beta\gamma\delta}(\xi,\eta, u,v,t)~=~\int_\R f\left(\xi,\eta,\tau-\mu\left(u \cdot\eta-v \cdot\xi\right)\right) |t-\tau|^{\big[{\alpha+\beta\over n+1}-{\gamma+\delta\over 2n+2}\big]-1} d\tau.
\eeq
Recall  the {\bf Hardy-Littlewood-Sobolev theorem} stated in the beginning of this paper. By applying (\ref{HLS Ineq}) with $\a={\alpha+\beta\over n+1}-{\gamma+\delta\over 2n+2}={1\over p}-{1\over q}$ and $\N=1$, we find
\bel{F regularity}
\begin{array}{lr}\ds
\left\{\int_{\R} \F^q_{\alpha\beta\gamma\delta}(\xi,\eta, u,v,t) dt \right\}^{1\over q}~\leq~\C_{p~q} \left\{\int_{\R} \Big[ f\left(\xi,\eta,t+\mu\left(u \cdot\eta-v \cdot\xi\right)\right)\Big]^p dt\right\}^{1\over p}
\\\\ \ds~~~~~~~~~~~~~~~~~~~~~~~~~~~~~~~~~~~~~~~~~~~~~
~=~\C_{p~q}~ \left\| f(\xi,\eta,\cdot)\right\|_{\L^p(\R)},\qquad (u,v)\in\R^n\times\R^n.
\end{array}
\eeq
From (\ref{I < =})-(\ref{F regularity}), we find
\bel{Est I-II}
\begin{array}{lr}\ds
\left\{\iiint_{\R^{2n+1}}{\sqrt{|u|^2+|v|^2}}^{-\gamma q} \Big(\I_{\alpha\beta\vartheta} f\Big)^q(u,v,t) dudvdt\right\}^{1\over q}
\\\\ \ds
\leq~\left\{ \iiint_{\R^{2n+1}}{\sqrt{|u|^2+|v|^2}}^{-\gamma q}\left\{\iint_{\R^{2n}} |u-\xi|^{\zeta-n}|v-\eta|^{\zeta-n} \F_{\alpha\beta\gamma\delta}(\xi,\eta, u,v,t)d\xi d\eta\right\}^q dudvdt\right\}^{1\over q}
\\\\ \ds
\leq~      \left\{\iint_{\R^{2n}}   {\sqrt{|u|^2+|v|^2}}^{-\gamma q}   \left\{ \iint_{\R^{2n}}   |u-\xi|^{\zeta-n}|v-\eta|^{\zeta-n}  \left\{\int_\R \F^q_{\alpha\beta\gamma\delta}(\xi,\eta, u,v,t) dt \right\}^{1\over q} d\xi d\eta\right\}^q dudv\right\}^{1\over q}
\\ \ds~~~~~~~~~~~~~~~~~~~~~~~~~~~~~~~~~~~~~~~~~~~~~~~~~~~~~~~~~~~~~~~~~~~~~~~~~~~~~~~~~~~~~~~~~~~~~~~~~~~
\hbox{\small{by Minkowski integral inequality}}
\\ \ds
\leq~\C_{p~q} \left\{\iint_{\R^{2n}}  {\sqrt{|u|^2+|v|^2}}^{-\gamma q}    \left\{ \iint_{\R^{2n}}   |u-\xi|^{\zeta-n}|v-\eta|^{\zeta-n}   \left\| f(\xi,\eta,\cdot)\right\|_{\L^p(\R)} d\xi d\eta\right\}^q dudv\right\}^{1\over q}.
\end{array}
\eeq
Define
\bel{II}
\II_\zeta g(u,v)~=~ \iint_{\R^{2n}} g(\xi,\eta)  |u-\xi|^{\zeta-n}|v-\eta|^{\zeta-n}    d\xi d\eta,\qquad 0<\zeta<n.
\eeq
Recall (\ref{Result Two})-(\ref{Formula Two}).
As a consequence of (\ref{Est I-II}), we can finish the proof of {\bf Theorem Two} by obtaining  the next two results.

{\bf Proposition One}~~{\it Let $\II_\zeta$ defined in (\ref{II}) for $0<\zeta<n$. Suppose $\omega(u,v)={\sqrt{|u|^2+|v|^2}}^{-\gamma}$, $\sigma(u,v)={\sqrt{|u|^2+|v|^2}}^\delta$ for $(u,v)\neq(0,0)$ and
$\gamma+\delta=0$.
We have
\bel{Prop One}
\begin{array}{cc}\ds
\left\| \omega \II_\zeta g\right\|_{\L^q(\R^{2n})}~\leq~\C_{p~q}~\left\| g\omega \right\|_{\L^p(\R^{2n})},\qquad 1<p< q<\infty
\\\\ \ds
\hbox{if}\qquad {\zeta\over n}~=~{1\over p}-{1\over q},\qquad \gamma<{n\over q},\qquad \delta<n\left({p-1\over p}\right).
\end{array}
\eeq
}

{\bf Proposition Two}~~{\it Let $\II_\zeta$ defined in (\ref{II}) for $0<\zeta<n$. Suppose $\omega(u,v)={\sqrt{|u|^2+|v|^2}}^{-\gamma}$, $\sigma(u,v)={\sqrt{|u|^2+|v|^2}}^\delta$ for $(u,v)\neq(0,0)$ and $\gamma+\delta>0$.
We have
\bel{Prop Two}
\begin{array}{cc}\ds
\left\| \omega \II_\zeta g\right\|_{\L^q(\R^{2n})}~\leq~\C_{p~q~\gamma~\delta}~\left\| g\sigma \right\|_{\L^p(\R^{2n})},\qquad 1<p\leq q<\infty
\\\\ \ds
\hbox{if}\qquad  \gamma<{2n\over q},\qquad \delta<2n\left({p-1\over p}\right),\qquad
{\zeta\over n}~=~{1\over p}-{1\over q}+{\gamma+\delta\over 2n};
\\\\ \ds
\zeta-{n\over p}<\delta ~~~\hbox{for}~~~ \gamma\ge0,~\delta\leq0;
\qquad
\zeta-n\left({q-1\over q}\right)<\gamma ~~~\hbox{for}~~~ \gamma\leq0,~\delta\ge0.
\end{array}
\eeq}

\subsection{Proof of Proposition One}
Observe that when $\gamma+\delta=0$, we have $\omega=\sigma$. Recall a classical one-weight theorem of fractional integrals due to Muckenhoupt and Wheeden \cite{Muckenhoupt-Wheeden}.

{\bf Muckenhoupt-Wheeden theorem} ~~{\it Let $\T_\a$ defined in (\ref{I_a}) for $0<\a<\N$. Suppose $\omega\ge0$ for $a.e ~x\in\R^\N$. Denote $\Q$ to be a cube in $\R^\N$ parallel to the coordinates.  We have
\bel{MW Ineq}
\left\| \omega \T_\a f\right\|_{\L^q(\R^\N)}~\leq~\C_{p~q}~\left\| f\omega\right\|_{\L^p(\R^\N)},\qquad 1<p<q<\infty
\eeq
if and only if
\bel{homo balance}
{\a\over \N}~=~{1\over p}-{1\over q}
\eeq
and
\bel{one-weight chara}
\left\{\vol\{\Q\}^{-1}\int_\Q \omega^q(x)dx\right\}^{1\over q} \left\{\vol\{\Q\}^{-1}\int_\Q \omega^{-{p\over p-1}}(x)dx\right\}^{p-1\over p}<\infty
\eeq 
for every $\Q\subset\R^\N$.}

Consider $\omega(u,v)={\sqrt{|u|^2+|v|^2}}^{-\gamma}$ where $\gamma+\delta=0$ for $\gamma<{n\over q}$ and $\delta<n\left({p-1\over p}\right)$. Take into account for $\a=\zeta$ and $\N=n$. For every $\Q\subset\R^n$, we simultaneously find
\bel{chara u}
\left\{\vol\{\Q\}^{-1}\int_\Q \Big[\sqrt{|u|^2+|v|^2}\Big]^{-\gamma q} du\right\}^{1\over q} \left\{\vol\{\Q\}^{-1}\int_\Q \Big[\sqrt{|u|^2+|v|^2}\Big]^{-\delta {p\over p-1}}du\right\}^{p-1\over p}<\infty,~~~~ v\in\R^n;
\eeq
\bel{chara v}
\left\{\vol\{\Q\}^{-1}\int_\Q \Big[\sqrt{|u|^2+|v|^2}\Big]^{-\gamma q} dv\right\}^{1\over q} \left\{\vol\{\Q\}^{-1}\int_\Q \Big[\sqrt{|u|^2+|v|^2}\Big]^{-\delta {p\over p-1}}dv\right\}^{p-1\over p}<\infty,~~~~ u\in\R^n.
\eeq
Indeed, by using $\gamma+\delta=0$, a standard one-parameter dilations in the left-hand -side of (\ref{chara u}) or (\ref{chara v}) shows that it is suffice to assume $\vol\{\Q\}^{1\over n}=1$. Moreover, $\omega^{q\gamma}(\cdot,v)$ and $\omega^{-\delta {p\over p-1}}(\cdot,v)$ are locally integrable in $\R^n$ for every $v\in\R^n$ provided that $\gamma<{n\over q}$ and $\delta<n\left({p-1\over p}\right)$. Vice versa for $\omega^{q\gamma}(u,\cdot)$ and $\omega^{-\delta {p\over p-1}}(u,\cdot)$. 
Let $\II_\zeta$ defined in (\ref{II}) for $0<\zeta<n$ and $g\ge0$. By applying {\bf Muckenhoupt-Wheeden theorem} two times, we have
\bel{EST one-weight}
\begin{array}{lr}\ds
\left\|\omega \II_\zeta g\right\|_{\L^q(\R^{2n})}~=~ \left\{\iint_{\R^{2n}}  {\sqrt{|u|^2+|v|^2}}^{-\gamma q}    \left\{ \iint_{\R^{2n}}   |u-\xi|^{\zeta-n}|v-\eta|^{\zeta-n}   g(\xi,\eta) d\xi d\eta\right\}^q dudv\right\}^{1\over q}
\\\\ \ds
~\leq~\C_{p~q} \left\{\int_{\R^n}       \left\{\int_{\R^n}  \left\{ \int_{\R^n}   |u-\xi|^{\zeta-n}  g(\xi,v) d\xi\right\}^p {\sqrt{|u|^2+|v|^2}}^{\delta p}dv\right\}^{q\over p} du\right\}^{1\over q}
\\\\ \ds
~\leq~\C_{p~q} \left\{\int_{\R^n} \left\{\int_{\R^n}  {\sqrt{|u|^2+|v|^2}}^{-\gamma q}       \left\{ \int_{\R^n}   |u-\xi|^{\zeta-n}  g(\xi,v)d\xi\right\}^q du \right\}^{p\over q} dv\right\}^{1\over p}
\\ \ds~~~~~~~~~~~~~~~~~~~~~~~~~~~~~~~~~~~~~~~~~~~~~~
\hbox{\small{by Minkowski integral inequality and $\gamma=-\delta$}}
\\ \ds
~\leq~\C_{p~q}      \left\{\iint_{\R^{2n}}  \Big[g(u,v)\Big]^p{\sqrt{|u|^2+|v|^2}}^{\delta p}du dv\right\}^{1\over p}~=~\C_{p~q}\left\| g\omega\right\|_{\L^p(\R^{2n})}.
\end{array}
\eeq

\section{Cone decomposition on $\R^n\times\R^n$}
\setcounter{equation}{0}
Let $\II_\zeta$ defined in (\ref{II}) for $0<\zeta<n$. For every $j\in\Z$, we consider
\bel{Delta_j II}
\Delta_j\II_\zeta g(u,v)~=~\iint_{\Lambda_j(u,v)} g(\xi,\eta) \left({1\over |u-\xi|}\right)^{n-\zeta}  \left({1\over |v-\eta|}\right)^{n-\zeta} d\xi d\eta
\eeq
where
\bel{Cone}
\Lambda_j(u,v)~=~\Bigg\{ (\xi,\eta)\in\R^n\times\R^n\colon 2^{-j}\leq{|u-\xi|\over |v-\eta|}<2^{-j+1}\Bigg\}.
\eeq
Observe that each $\Lambda_j(u,v)$ is a dyadic cone centered on $(u,v)\in\R^n\times\R^n$ with an eccentricity depending on $j\in\Z$. 
\begin{figure}[h]
\centering
\includegraphics[scale=0.40]{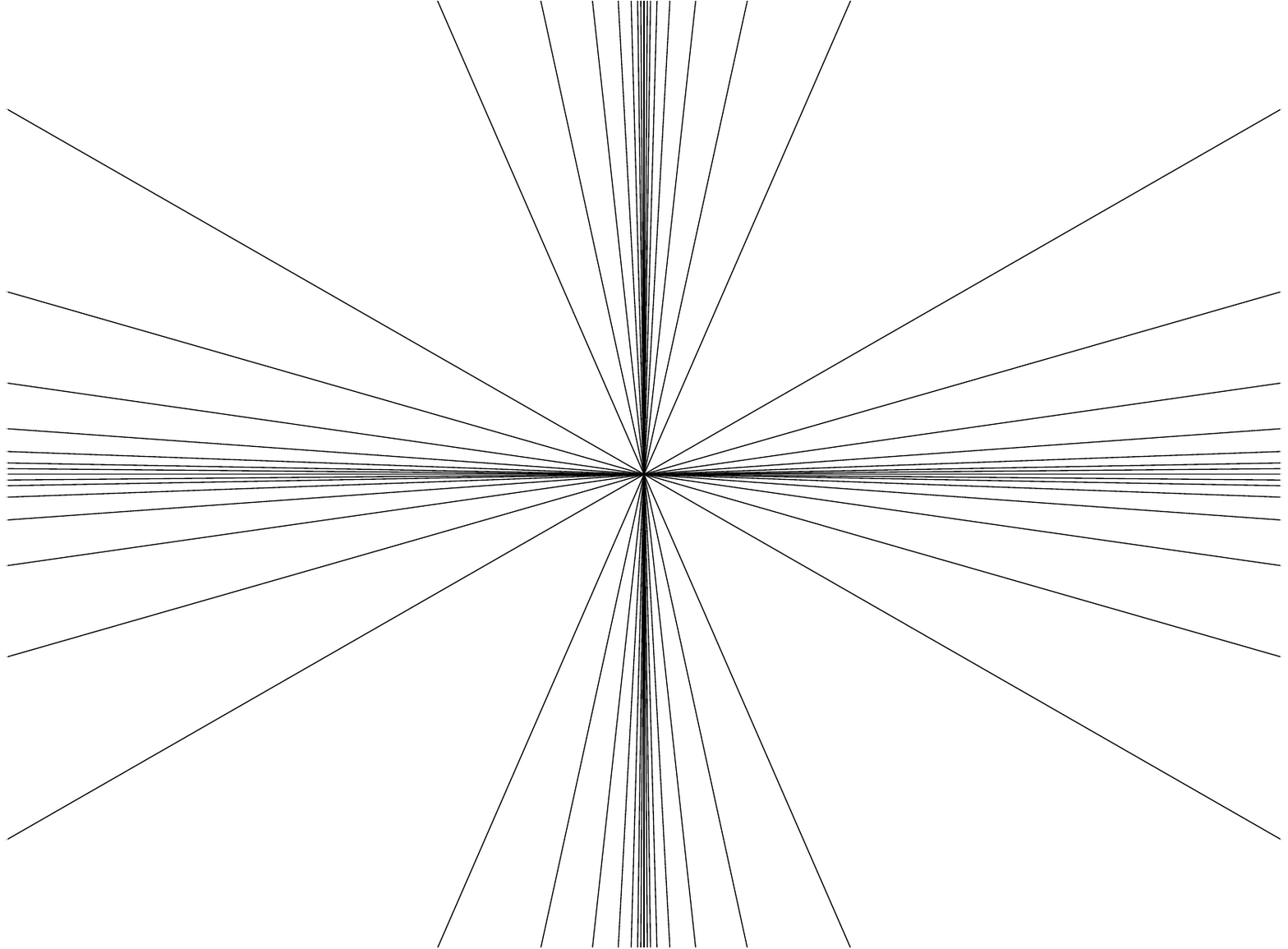}
\end{figure}

Denote $\Q^j_i$ to be a dilated of $\Q_i\subset\R^n$ such that
$\vol\{\Q_i^j\}^{1\over n}=2^{-j} \vol\{\Q_i\}^{1\over n}$ for $i=1,2$ and $j\in\Z$. Let $r>1$.
We have
\bel{A-Characteristic Dila}
\begin{array}{lr}\ds
\prod_{i=1}^2 \vol\{\Q_i\}^{{\zeta\over n}-{1\over p}+{1\over q}}\left\{{1\over\vol\{\Q_1\}\vol\{\Q_2\}}\iint_{\Q_1\times\Q_2} \omega^{qr}\left(2^{-j} u,v\right)dudv\right\}^{1\over qr}
\\ \ds~~~~~~~~~~~~~~~~~~~~~~~~~~~~
\left\{{1\over\vol\{\Q_1\}\vol\{\Q_2\}}\iint_{\Q_1\times\Q_2} \left({1\over \sigma}\right)^{pr\over p-1}\left(2^{-j} u,v\right)dudv\right\}^{p-1\over pr}
\\\\ \ds
=~2^{j\big[\zeta-{n\over p}+{n\over q}\big]} \vol\{\Q_1^j\}^{{\zeta\over n}-{1\over p}+{1\over q}}
\vol\{\Q_2\}^{{\zeta\over n}-{1\over p}+{1\over q}}
\\ \ds
\left\{{1\over\vol\{\Q_1^j\}\vol\{\Q_2\}}\iint_{\Q_1^j\times\Q_2} \omega^{qr}\left( u,v\right)dudv\right\}^{1\over qr}
\left\{{1\over\vol\{\Q_1^j\}\vol\{\Q_2\}}\iint_{\Q_1^j\times\Q_2} \left({1\over \sigma}\right)^{pr\over p-1}\left( u,v\right)dudv\right\}^{p-1\over pr}.
\end{array}
\eeq

Given $j\in\Z$, we define
\bel{sup A_pqr^alpha t}
\begin{array}{lr}\ds
\A_{pqr}^\zeta\left(j~\colon\omega,\sigma\right)~=~\sup_{\Q_1\times\Q_2\subset\R^n\times\R^n\colon~ \vol\{\Q_1\}^{1\over n}/\vol\{\Q_2\}^{1\over n}=2^{-j}} \prod_{i=1}^2 \vol\{\Q_i\}^{{\zeta\over n}-{1\over p}+{1\over q}}
\\ \ds
\left\{{1\over\vol\{\Q_1\}\vol\{\Q_2\}}\iint_{\Q_1\times\Q_2} \omega^{qr}\left( u,v\right)dudv\right\}^{1\over qr}
\left\{{1\over\vol\{\Q_1\}\vol\{\Q_2\}}\iint_{\Q_1\times\Q_2} \left({1\over \sigma}\right)^{pr\over p-1}\left( u,v\right)dudv\right\}^{p-1\over pr}.
\end{array}
\eeq
Suppose $\vol\{\Q_1\}^{1\over n}=\vol\{\Q_2\}^{1\over n}$.
We find 
\bel{A-Characteristic Bound}
\begin{array}{lr}\ds
\prod_{i=1}^2 \vol\{\Q_i\}^{{\zeta\over n}-{1\over p}+{1\over q}}\left\{{1\over\vol\{\Q_1\}\vol\{\Q_2\}}\iint_{\Q_1\times\Q_2} \omega^{qr}\left(2^{-j} u,v\right)dudv\right\}^{1\over qr}
\\ \ds~~~~~~~~~~~~~~~~~~~~~~~~~~~~
\left\{{1\over\vol\{\Q_1\}\vol\{\Q_2\}}\iint_{\Q_1\times\Q_2} \left({1\over \sigma}\right)^{pr\over p-1}\left(2^{-j} u,v\right)dudv\right\}^{p-1\over pr}
\\\\ \ds
~\leq~2^{j\big[\zeta-{n\over p}+{n\over q}\big]} \A_{pqr}^\zeta\left(j~\colon\omega,\sigma\right)\qquad\hbox{\small{by (\ref{A-Characteristic Dila})-(\ref{sup A_pqr^alpha t})}}.
\end{array}
\eeq
Next,  recall a classical result due to Sawyer and Wheeden \cite{Sawyer-Wheeden}  for  one-parameter fractional integral operators in weighted norms. 

Suppose 
\bel{A0 sup}
\begin{array}{lr}\ds
\A_{pqr}^\zeta(0\colon \omega, \sigma)~<~\infty\qquad\hbox{\small{for some $r>1$}}.
\end{array}
\eeq
We have
\bel{One-para Norm Ineq}
\begin{array}{lr}\ds
\left\{\iint_{\R^{2n}}\left\{\iint_{\R^{2n}}g( \xi,\eta)\left[{1\over \sqrt{|u-\xi|^2+|v-\eta|^2}}\right]^{2n-2\zeta} d\xi d\eta\right\}^q \omega^q(u,v)dudv\right\}^{1\over q} 
\\\\ \ds
~\leq~\C_{p~q~r~\zeta}~ \A_{pqr}^\zeta(0\colon \omega, \sigma)~\left\{\iint_{\R^{2n}} \Big(g\sigma\Big)^p(u,v)dudv\right\}^{1\over p},\qquad 1<p\leq q<\infty.
\end{array}
\eeq
 \begin{remark} 
The  constant $\C_{p~q~r~\zeta} ~\A_{pqr}^\zeta(0\colon \omega, \sigma)$ in (\ref{One-para Norm Ineq}) is not written explicitly in the original statement by Saywer and Wheeden \cite{Sawyer-Wheeden} (Theorem 1). But, it can be computed directly by carrying out the proof given in section 2 of \cite{Sawyer-Wheeden}.\end{remark}
By applying (\ref{One-para Norm Ineq})  and  using the estimate in (\ref{A-Characteristic Bound}), we find
\bel{Regularity est}
\begin{array}{lr}\ds
\left\{\iint_{\R^{2n}}\left\{\iint_{\R^{2n}}g(2^{-j}\xi,\eta)\left[{1\over \sqrt{|u-\xi|^2+|v-\eta|^2}}\right]^{2n-2\zeta} d\xi d\eta\right\}^{q}\omega^q(2^{-j} u,v)dudv\right\}^{1\over q} 
\\\\ \ds
~\leq~\C_{p~q~r~\zeta}~2^{j\big[\zeta-{n\over p}+{n\over q}\big]} \A_{pqr}^\zeta\left(j~\colon\omega,\sigma\right)\left\{\iint_{\R^{2n}} \Big(g\sigma\Big)^p(2^{-j}u,v)dudv\right\}^{1\over p}
\end{array}
\eeq
for $1<p\leq q<\infty$ and every $j\in\Z$.

Recall (\ref{Delta_j II})-(\ref{Cone}). By changing dilations $(u,v)\mt (2^{-j}u,v)$ and $(\xi,\eta)\mt (2^{-j}\xi,\eta)$, we have
\bel{Dilation Est}
\begin{array}{lr}\ds
\left\{\iint_{\R^{2n}} \Big(\Delta_j\II_\zeta g\Big)^q(u,v)\omega^q(u,v)dudv\right\}^{1\over q}
\\\\ \ds
~=~\left\{\iint_{\R^{2n}}\left\{ \iint_{\Lambda_j(u,v)}g(\xi,\eta)\left({1\over |u-\xi|}\right)^{n-\zeta}\left({1\over |v-\eta|}\right)^{n-\zeta}d\xi d\eta\right\}^q\omega^q(u,v)dudv\right\}^{1\over q}
 \\\\ \ds
 ~=~\left\{\iint_{\R^{2n}}\left\{\iint_{\Lambda_0(u,v)}g(2^{-j}\xi,\eta)\left({1\over 2^{-j}|u-\xi|}\right)^{n-\zeta} \left({1\over |v-\eta|}\right)^{n-\zeta}2^{-jn} d\xi d\eta\right\}^q\omega^q(2^{-j} u,v) 2^{-jn} dudv\right\}^{1\over q}
 \\\\ \ds
 ~\lesssim~~2^{-j\big[\zeta+{n\over q}\big]}\left\{\iint_{\R^{2n}}\left\{ \iint_{\R^{2n}}g(2^{-j} \xi,\eta)\left[{1\over \sqrt{|u-\xi|^2+|v-\eta|^2}}\right]^{2n-2\zeta}d\xi d\eta\right\}^q\omega^q(2^{-j} u,v)dudv\right\}^{1\over q}
 \\\\ \ds
 ~\leq~\C_{p~q~r~\zeta}~2^{-j\big[\zeta+{n\over q}\big]}2^{j\big[\zeta-{n\over p}+{n\over q}\big]}\A_{pqr}^\zeta\left(j~\colon\omega,\sigma\right)\left\{ \iint_{\R^{2n}} \Big(g\sigma\Big)^p(2^{-j} u,v)dudv\right\}^{1\over p}\qquad \hbox{\small{by (\ref{Regularity est})}}
  \\\\ \ds
   ~=~\C_{p~q~r~\zeta}~\A_{pqr}^\zeta(j~\colon\omega,\sigma)~ 2^{-j\big[\zeta+{n\over q}\big]}2^{j\big[\zeta-{n\over p}+{n\over q}\big]}\left\{ \iint_{\R^{2n}} \Big(g\sigma\Big)^p\left(u,v\right) 2^{jn}dudv\right\}^{1\over p}
 \\\\ \ds
 ~=~\C_{p~q~r~\zeta}~\A_{pqr}^\zeta\left(j~\colon\omega,\sigma\right)\left\{ \iint_{\R^{2n}} \Big(g\sigma\Big)^p(u,v)dxudv\right\}^{1\over p}.
\end{array}
\eeq
By using (\ref{Dilation Est}) and Minkowski inequality, we obtain the $\L^p\mt\L^q$-norm inequality in (\ref{Result Two}) provided that
\[
\sum_{j\in\Z} \A_{pqr}^\zeta(j~\colon\omega,\sigma)~<~\infty.
\]

{\bf Principal Lemma}~~{\it Suppose $\omega(u,v)={\sqrt{|u|^2+|v|^2}}^{-\gamma}$ and $\sigma(u,v)={\sqrt{|u|^2+|v|^2}}^\delta$ for $(u,v)\neq(0,0)$. Let $\gamma+\delta>0$. 
There exists  $\ve=\ve(p,q,\gamma,\delta)>0$  such that
 \bel{Summability}
 \begin{array}{cc}\ds
 \A_{pqr}^\zeta(j~\colon\omega,\sigma)~<~\C_{p~q~\gamma~\delta} ~2^{-\ve |j|}
 \\\\ \ds
\hbox{if}\qquad  \gamma<{2n\over q},\qquad \delta<2n\left({p-1\over p}\right),\qquad
{\zeta\over n}~=~{1\over p}-{1\over q}+{\gamma+\delta\over 2n};
\\\\ \ds
\zeta-{n\over p}<\delta \qquad\hbox{for}\qquad \gamma\ge0,~\delta\leq0;
\\\\ \ds
\zeta-n\left({q-1\over q}\right)<\gamma \qquad\hbox{for}\qquad \gamma\leq0,~\delta\ge0 
 \end{array}
\eeq
for some $r=r(p,q,\gamma,\delta)>1$ and every $j\in\Z$.}

By symmetry, we consider $j>0$ only. For every $\Q_1\times\Q_2\subset\R^n\times\R^n$ satisfying 
\bel{ratio Q}
\vol\{\Q_1\}^{1\over n}/\vol\{\Q_2\}^{1\over n}~=~\lambda,\qquad 0<\lambda\leq1,
\eeq
we aim to show that the constraints of $p,q,\gamma,\delta$ inside (\ref{Summability}) imply
\bel{Eccentricity Decay}
\begin{array}{lr}\ds
\prod_{i=1}^2 \vol\{\Q_i\}^{{\zeta\over n}-{1\over p}+{1\over q}}\left\{{1\over\vol\{\Q_1\}\vol\{\Q_2\}}\iint_{\Q_1\times\Q_2} \left({1\over|u|+|v|}\right)^{\gamma qr}dudv\right\}^{1\over qr}
\\ \ds~~~~~~~~~~~~~~~~~~~~~~~~~~~~
\left\{{1\over\vol\{\Q_1\}\vol\{\Q_2\}}\iint_{\Q_1\times\Q_2} \left({1\over |u|+|v|}\right)^{\delta{pr\over p-1}}dudv\right\}^{p-1\over pr}
~\leq~ \C_{p~q~r~\gamma~\delta}~\lambda^\ve
\end{array}
\eeq
where $\ve>0$ and $r>1$  depend on $p,q,\gamma, \delta$. 

By  using the homogeneity condition ${\zeta\over n}={1\over p}-{1\over q}+{\gamma+\delta\over 2n}$, we find that the left-hand-side of  (\ref{Eccentricity Decay}) is invariant by changing  dilations in one-parameter. Therefore, it is  suffice to assert $\vol\{\Q_2\}^{1\over n}=1$.
\begin{remark}  Let $\Q_i^o$ and $\Q_i^*\subset\R^n$  be cubes centered on the origin of $\R^n$ and 
\bel{Q*}
\vol\{\Q_i^o\}^{1\over n}~=~\vol\{\Q_i\}^{1\over n},\qquad \vol\{\Q_i^*\}^{1\over n}~=~3\vol\{\Q_i\}^{1\over n},\qquad i~=~1,2.
\eeq
Suppose $\Q_i\cap\Q_i^o=\emptyset$. We must have $|x|\ge|x^o|/\sqrt{n}$ for every $x\in\Q_i$ and $x^o\in\Q_i^o$. 

Otherwise,   $\Q_i\subset\Q_i^*$ if $\Q_i$ intersects $\Q_i^o$.
\end{remark}
Suppose $\Q_1\times\Q_2$ centered on $(u_o,v_o)\in\R^n\times\R^n$ of which $\sqrt{|u_o|^2+|v_o|^2}>3$. Because $\Q_1\times\Q_2$ has a diameter $1$, we find
\[{1\over 2} \sqrt{|u_o|^2+|v_o|^2}~\leq~\sqrt{|u|^2+|v|^2}~\leq~2\sqrt{|u_o|^2+|v_o|^2},\qquad (u,v)\in\Q_1\times\Q_2.\]
This further implies
\bel{Eccentricity decay >}
\begin{array}{lr}\ds
\prod_{i=1}^2 \vol\{\Q_i\}^{{\zeta\over n}-{1\over p}+{1\over q}}\left\{{1\over\vol\{\Q_1\}\vol\{\Q_2\}}\iint_{\Q_1\times\Q_2} \left({1\over|u|+|v|}\right)^{\gamma qr}dudv\right\}^{1\over qr}
\\ \ds~~~~~~~~~~~~~~~~~~~~~~~~~~~~
\left\{{1\over\vol\{\Q_1\}\vol\{\Q_2\}}\iint_{\Q_1\times\Q_2} \left({1\over |u|+|v|}\right)^{\delta{pr\over p-1}}dudv\right\}^{p-1\over pr}
\\\\ \ds
~\leq~ \C_{p~q~r~\gamma~\delta}~\Big[\sqrt{|u_o|^2+|v_o|^2}\Big]^{-(\gamma+\delta)} \lambda^{{\zeta\over n}-{1\over p}+{1\over q}}
\\\\ \ds
~\leq~ \C_{p~q~r~\gamma~\delta}~ \lambda^\ve,\qquad \hbox{\small{$\ve={\zeta\over n}-{1\over p}+{1\over q}={\gamma+\delta\over 2n}>0$}}.
\end{array}
\eeq
\begin{remark}
From now on, we assume $\Q_1\times\Q_2$ centered on some $(u_o,v_o)\in\R^n\times\R^n$ with $\sqrt{|u_o|^2+|v_o|^2}\leq3$. 
\end{remark}
Let $\Q^*_1$ defined (\ref{Q*}). We have
 \bel{Int Compara gamma}
 \begin{array}{lr}\ds
 \int_{\Q_1} \left({1\over |u|}\right)^{\gamma qr}  du~\lesssim~ \int_{\Q_1^*} \left({1\over |u|}\right)^{\gamma qr}  du,\qquad \hbox{$0<\gamma qr<n$};
  \\\\ \ds
  \int_{\Q_1} \left({1\over |u|}\right)^{\gamma qr-n}  du~\lesssim~ \int_{\Q_1^*} \left({1\over |u|}\right)^{\gamma qr-n}  du,\qquad  \hbox{$n<\gamma qr<2n$}
   \end{array}
 \eeq
 and
  \bel{Int Compara delta}
 \begin{array}{lr}\ds
 \int_{ \Q_1}  \left({1\over |u|}\right)^{\delta{pr\over p-1}} du~\lesssim~\int_{ \Q_1^*}  \left({1\over |u|}\right)^{\delta{pr\over p-1}} du ,\qquad \hbox{$0<\delta \left({p\over p-1}\right)r<n$};
 \\\\ \ds
 \int_{ \Q_1}  \left({1\over |u|}\right)^{\delta{pr\over p-1}-n} du~\lesssim~\int_{ \Q_1^*}  \left({1\over |u|}\right)^{\delta{pr\over p-1}-n} du, \qquad \hbox{$n<\delta \left({p\over p-1}\right)r<2n$}
   \end{array}
 \eeq

The remaining proof  is split into 3 cases, $w.r.t$ $\gamma\ge0,\delta\leq0$; $\gamma\leq0,\delta\ge0$ and $\gamma>0,\delta>0$.

\subsection{Case One: $\gamma\ge0$, $\delta\leq0$}
By adjusting the value of $r>1$, we find
\bel{Range 1}
\begin{array}{cc}\ds
\hbox{$0<\gamma qr<n$}\qquad \hbox{or}\qquad \hbox{$n<\gamma qr<2n$}.
\end{array}
 \eeq
Suppose $0<\gamma qr<n$. We have
\bel{Decay Case1 Est1}
\begin{array}{lr}\ds
\prod_{i=1}^2 \vol\{\Q_i\}^{{\zeta\over n}-{1\over p}+{1\over q}}\left\{{1\over\vol\{\Q_1\}\vol\{\Q_2\}}\iint_{\Q_1\times\Q_2} \left({1\over|u|+|v|}\right)^{\gamma qr}dudv\right\}^{1\over qr}
\\ \ds~~~~~~~~~~~~~~~~~~~~~~~~~~~~
\left\{{1\over\vol\{\Q_1\}\vol\{\Q_2\}}\iint_{\Q_1\times\Q_2} \left({1\over |u|+|v|}\right)^{\delta{pr\over p-1}}dudv\right\}^{p-1\over pr}
\\\\ \ds
\leq~\C_{p~r~\delta}~ \lambda^{\zeta-{n\over p}+{n\over q}}\left({1\over \lambda}\right)^{n\over qr}  
\left\{\int_{\Q_1}\left\{\int_{\Q_2}  \left({1\over |u|+|v|}\right)^{\gamma qr} dv\right\} du\right\}^{1\over qr}\qquad\hbox{\small{ by {\bf Remark 5.3}  ($\delta\leq0$)}}

\\\\ \ds
\leq~ \C_{p~r~\delta}~\lambda^{\zeta-{n\over p}+{n\over q}}\left({1\over \lambda}\right)^{n\over qr}  \lambda^{n\over qr} 
\left\{\int_{\Q_2}  \left({1\over |v|}\right)^{\gamma qr} dv\right\}^{1\over qr}

\\\\ \ds
\leq~\C_{p~r~\delta}~ \lambda^{\zeta-{n\over p}+{n\over q}}\left({1\over \lambda}\right)^{n\over qr}   \lambda^{n\over qr} 
\left\{\int_{\Q_2^*}  \left({1\over |v|}\right)^{\gamma qr} dv\right\}^{1\over qr}
\qquad \hbox{\small{by (\ref{Int Compara gamma})}}
\\\\ \ds
\leq~\C_{p~q~r~\gamma~\delta} ~\lambda^{\zeta-{n\over p}+{n\over q}}\left({1\over \lambda}\right)^{n\over qr}  
\lambda^{n\over qr}

\\\\ \ds
=~\C_{p~q~r~\gamma~\delta}~ \lambda^{\zeta-{n\over p}+{n\over q}}
\\\\ \ds
=~\C_{p~q~r~\gamma~\delta}~\lambda^{\gamma+\delta\over 2} \qquad\hbox{( \small{${\zeta\over n}={1\over p}-{1\over q}+{\gamma+\delta\over 2n}$} )}
\\\\ \ds
=~\C_{p~q~r~\gamma~\delta}~\lambda^\ve,\qquad \hbox{\small{$\ve={\gamma+\delta\over 2}>0$}}.
\end{array}
\eeq  
 Suppose $n<\gamma qr<2n$. Recall  $\zeta-{n\over p}<\delta$ as an necessity. Together with the homogeneity condition ${\zeta\over n}={1\over p}-{1\over q}+{\gamma+\delta\over 2n}$, we find
 \bel{Case One compu}
 \begin{array}{lr}\ds
 \zeta-{n\over p}~=~-{n\over q}+{\gamma+\delta\over 2}~<~\delta
 \\\\ \ds
 \Longrightarrow\qquad {n\over q}-{\gamma\over 2}+{\delta\over 2}~>~0.
 \end{array}
 \eeq
For $r$ chosen sufficiently close to $1$, we have
\bel{Decay Case1 Est2}
\begin{array}{lr}\ds
\prod_{i=1}^2 \vol\{\Q_i\}^{{\zeta\over n}-{1\over p}+{1\over q}}\left\{{1\over\vol\{\Q_1\}\vol\{\Q_2\}}\iint_{\Q_1\times\Q_2} \left({1\over|u|+|v|}\right)^{\gamma qr}dudv\right\}^{1\over qr}
\\ \ds~~~~~~~~~~~~~~~~~~~~~~~~~~~~
\left\{{1\over\vol\{\Q_1\}\vol\{\Q_2\}}\iint_{\Q_1\times\Q_2} \left({1\over |u|+|v|}\right)^{\delta{pr\over p-1}}dudv\right\}^{p-1\over pr}
\\\\ \ds
\leq~ \C_{p~r~\delta}~\lambda^{\zeta-{n\over p}+{n\over q}}\left({1\over \lambda}\right)^{n\over qr}  
\left\{\int_{\Q_1}\left\{\int_{\Q_2}  \left({1\over |u|+|v|}\right)^{\gamma qr} dv\right\} du\right\}^{1\over qr}
\qquad\hbox{\small{by {\bf Remark 5.3} ($\delta\leq0$)}}
\\\\ \ds
\leq~\C_{p~r~\delta}~ \lambda^{\zeta-{n\over p}+{n\over q}}\left({1\over \lambda}\right)^{n\over qr}  
\left\{\int_{\Q_1}\left\{\int_{\R^n}  \left({1\over |u|+|v|}\right)^{\gamma qr} dv\right\} du\right\}^{1\over qr}
\\\\ \ds
\approx~  \C_{p~r~\delta}~\lambda^{\zeta-{n\over p}+{n\over q}}\left({1\over \lambda}\right)^{n\over qr}  
\\ \ds~~~
\left\{\int_{\Q_1}\left\{\idotsint_{\R^n}  \left({1\over |u|+|v_1|+\cdots+|v_n|}\right)^{\gamma qr} dv_1\cdots dv_n\right\} du\right\}^{1\over qr}
\\\\ \ds
\leq~\C_{p~q~r~\gamma~\delta}~ \lambda^{\zeta-{n\over p}+{n\over q}}\left({1\over \lambda}\right)^{n\over qr}  
\left\{\int_{\Q_1} \left({1\over |u|}\right)^{\gamma qr-n}  du\right\}^{1\over qr}
\\\\ \ds
\leq~\C_{p~q~r~\gamma~\delta} ~\lambda^{\zeta-{n\over p}+{n\over q}}\left({1\over \lambda}\right)^{n\over qr}  
\left\{\int_{\Q_1^*} \left({1\over |u|}\right)^{\gamma qr-n}  du\right\}^{1\over qr}
\qquad
 \hbox{\small{by (\ref{Int Compara gamma})}}

\\\\ \ds
\leq~\C_{p~q~r~\gamma~\delta}~ \lambda^{\zeta-{n\over p}+{n\over q}}\left({1\over \lambda}\right)^{n\over qr}      \lambda^{{2n\over qr}-\gamma}
\\\\ \ds
=~\C_{p~q~r~\gamma~\delta}~\lambda^{\zeta-{n\over p}+{n\over q}}\lambda^{{n\over qr}-\gamma}
~=~\C_{p~q~r~\gamma~\delta}~\lambda^{\gamma+\delta\over 2}\lambda^{{n\over qr}-\gamma} \qquad\hbox{( \small{${\zeta\over n}={1\over p}-{1\over q}+{\gamma+\delta\over 2n}$} )}
\\\\ \ds
=~\C_{p~q~r~\gamma~\delta}~\lambda^{{n\over qr}-{\gamma\over 2}+{\delta\over 2}}
\\\\ \ds
=~\C_{p~q~r~\gamma~\delta}~\lambda^\ve,\qquad \hbox{\small{$\ve={n\over qr}-{\gamma\over 2}+{\delta\over 2}>0$ ~~by (\ref{Case One compu})}}.
\end{array}
\eeq

 \subsection{Case Two: $\gamma\leq0$, $\delta\ge0$}
By adjusting the value of $r>1$, we find
\bel{Range 2}
\begin{array}{cc}\ds
\hbox{$0<\delta \left({p\over p-1}\right)r<n$}\qquad \hbox{or}\qquad \hbox{$n<\delta \left({p\over p-1}\right)r<2n$}.
\end{array}
 \eeq
Suppose $0<\delta \left({p\over p-1}\right)r<n$. We have
\bel{Decay Case2 Est1}
\begin{array}{lr}\ds
\prod_{i=1}^2 \vol\{\Q_i\}^{{\zeta\over n}-{1\over p}+{1\over q}}\left\{{1\over\vol\{\Q_1\}\vol\{\Q_2\}}\iint_{\Q_1\times\Q_2} \left({1\over|u|+|v|}\right)^{\gamma qr}dudv\right\}^{1\over qr}
\\ \ds~~~~~~~~~~~~~~~~~~~~~~~~~~~~
\left\{{1\over\vol\{\Q_1\}\vol\{\Q_2\}}\iint_{\Q_1\times\Q_2} \left({1\over |u|+|v|}\right)^{\delta{pr\over p-1}}dudv\right\}^{p-1\over pr}
\\\\ \ds
\leq~\C_{q~r~\gamma}~ \lambda^{\zeta-{n\over p}+{n\over q}}\left({1\over \lambda}\right)^{n\left({p-1\over pr}\right)}  
\left\{\int_{\Q_1}\left\{\int_{\Q_2}  \left({1\over |u|+|v|}\right)^{\delta {pr\over p-1}} dv\right\} du\right\}^{p-1\over pr}\qquad\hbox{\small{by {\bf Remark 5.3} ($\gamma\leq0$)}}

\\\\ \ds
\leq~ \C_{q~r~\gamma}~\lambda^{\zeta-{n\over p}+{n\over q}}\left({1\over \lambda}\right)^{n\left({p-1\over pr}\right)}  \lambda^{n\left({p-1\over pr}\right)} 
\left\{\int_{\Q_2}  \left({1\over |v|}\right)^{\delta {pr\over p-1}} dv\right\}^{p-1\over pr}

\\\\ \ds
\leq~ \C_{q~r~\gamma}~\lambda^{\zeta-{n\over p}+{n\over q}}\left({1\over \lambda}\right)^{n\left({p-1\over pr}\right)}  \lambda^{n\left({p-1\over pr}\right)} 
\left\{\int_{\Q_2^*}  \left({1\over |v|}\right)^{\delta {pr\over p-1}} dv\right\}^{p-1\over pr}
\qquad \hbox{\small{by (\ref{Int Compara delta})}}
\\\\ \ds
\leq~\C_{p~q~r~\gamma~\delta} ~\lambda^{\zeta-{n\over p}+{n\over q}}\left({1\over \lambda}\right)^{n\left({p-1\over pr}\right)}  \lambda^{n\left({p-1\over pr}\right)} 

\\\\ \ds
=~\C_{p~q~r~\gamma~\delta}~ \lambda^{\zeta-{n\over p}+{n\over q}}
\\\\ \ds
=~\C_{p~q~r~\gamma~\delta}~\lambda^{\gamma+\delta\over 2} \qquad\hbox{( \small{${\zeta\over n}={1\over p}-{1\over q}+{\gamma+\delta\over 2n}$} )}
\\\\ \ds
=~\C_{p~q~r~\gamma~\delta}~\lambda^\ve,\qquad \hbox{\small{$\ve={\gamma+\delta\over 2}>0$}}.
\end{array}
\eeq  
 Suppose $n<\delta \left({p\over p-1}\right)r<2n$. Recall  $\zeta-n\left({q-1\over q}\right)<\gamma$ as an necessity. Together with the homogeneity condition ${\zeta\over n}={1\over p}-{1\over q}+{\gamma+\delta\over 2n}={q-1\over q}-{p-1\over p}+{\gamma+\delta\over 2n}$, we find
 \bel{Case Two compu}
 \begin{array}{lr}\ds
 \zeta-n\left({q-1\over q}\right)~=~-n\left({p-1\over p}\right)+{\gamma+\delta\over 2}~<~\gamma
 \\\\ \ds
 \Longrightarrow\qquad n\left({p-1\over p}\right)+{\gamma\over 2}-{\delta\over 2}~>~0.
 \end{array}
 \eeq
For $r$ chosen sufficiently close to $1$, we have
\bel{Decay Case2 Est2}
\begin{array}{lr}\ds
\prod_{i=1}^2 \vol\{\Q_i\}^{{\zeta\over n}-{1\over p}+{1\over q}}\left\{{1\over\vol\{\Q_1\}\vol\{\Q_2\}}\iint_{\Q_1\times\Q_2} \left({1\over|u|+|v|}\right)^{\gamma qr}dudv\right\}^{1\over qr}
\\ \ds~~~~~~~~~~~~~~~~~~~~~~~~~~~~
\left\{{1\over\vol\{\Q_1\}\vol\{\Q_2\}}\iint_{\Q_1\times\Q_2} \left({1\over |u|+|v|}\right)^{\delta{pr\over p-1}}dudv\right\}^{p-1\over pr}
\\\\ \ds
\leq~ \C_{q~r~\gamma}~\lambda^{\zeta-{n\over p}+{n\over q}}\left({1\over \lambda}\right)^{n\left({p-1\over pr}\right)}  
\left\{\int_{\Q_1}\left\{\int_{\Q_2}  \left({1\over |u|+|v|}\right)^{\delta {pr\over p-1}} dv\right\} du\right\}^{p-1\over pr}
\qquad\hbox{\small{by {\bf Remark 5.3} ($\gamma\leq0$)}}
\\\\ \ds
\leq~\C_{q~r~\gamma}~\lambda^{\zeta-{n\over p}+{n\over q}}\left({1\over \lambda}\right)^{n\left({p-1\over pr}\right)}  
\left\{\int_{\Q_1}\left\{\int_{\R^n}  \left({1\over |u|+|v|}\right)^{\delta {pr\over p-1}} dv\right\} du\right\}^{p-1\over pr}
\\\\ \ds
\approx~  \C_{q~r~\gamma}~\lambda^{\zeta-{n\over p}+{n\over q}}\left({1\over \lambda}\right)^{n\left({p-1\over pr}\right)}  
\\ \ds~~~
\left\{\int_{\Q_1}\left\{\idotsint_{\R^n}  \left({1\over |u|+|v_1|+\cdots+|v_n|}\right)^{\delta{pr\over p-1}} dv_1\cdots dv_n\right\} du\right\}^{p-1\over pr}
\\\\ \ds
\leq~\C_{p~q~r~\gamma~\delta}~ \lambda^{\zeta-{n\over p}+{n\over q}}\left({1\over \lambda}\right)^{n\left({p-1\over pr}\right)}  
\left\{\int_{\Q_1} \left({1\over |u|}\right)^{\delta {pr\over p-1}-n}  du\right\}^{p-1\over pr}
\\\\ \ds
\leq~\C_{p~q~r~\gamma~\delta}~ \lambda^{\zeta-{n\over p}+{n\over q}}\left({1\over \lambda}\right)^{n\left({p-1\over pr}\right)}  
\left\{\int_{\Q_1^*} \left({1\over |u|}\right)^{\delta {pr\over p-1}-n}  du\right\}^{p-1\over pr}
\qquad
 \hbox{\small{by (\ref{Int Compara delta})}}

\\\\ \ds
\leq~\C_{p~q~r~\gamma~\delta}~ \lambda^{\zeta-{n\over p}+{n\over q}}\left({1\over \lambda}\right)^{n\left({p-1\over pr}\right)}      \lambda^{2n\left({p-1\over pr}\right)-\delta}
\\\\ \ds
=~\C_{p~q~r~\gamma~\delta}~\lambda^{\zeta-{n\over p}+{n\over q}}\lambda^{n\left({p-1\over pr}\right)-\delta}
\\\\ \ds
=~\C_{p~q~r~\gamma~\delta}~\lambda^{\gamma+\delta\over 2}\lambda^{n\left({p-1\over pr}\right)-\delta} \qquad\hbox{( \small{${\zeta\over n}={1\over p}-{1\over q}+{\gamma+\delta\over 2n}$} )}
\\\\ \ds
=~\C_{p~q~r~\gamma~\delta}~\lambda^{n\left({p-1\over pr}\right)+{\gamma\over 2}-{\delta\over 2}}
\\\\ \ds
=~\C_{p~q~r~\gamma~\delta}~\lambda^\ve,\qquad \hbox{\small{$\ve=n\left({p-1\over pr}\right)+{\gamma\over 2}-{\delta\over 2}>0$ ~~by (\ref{Case Two compu})}}.
\end{array}
\eeq

\subsection{Case Three: $\gamma>0$, $\delta>0$}
By adjusting the value of $r>1$, we find
\bel{Range 3}
\begin{array}{cc}\ds
\hbox{$0<\gamma qr<n$},\qquad \hbox{$0<\delta \left({p\over p-1}\right)r<n$};
\\\\ \ds
\hbox{$n<\gamma qr<2n$},\qquad \hbox{$0<\delta \left({p\over p-1}\right)r<n$}\qquad\hbox{\small{or}}\qquad \hbox{$0<\gamma qr<n$},\qquad \hbox{$n<\delta \left({p\over p-1}\right)r<2n$};
\\\\ \ds
\hbox{$n<\gamma qr<2n$},\qquad \hbox{$n<\delta \left({p\over p-1}\right)r<2n$}.
\end{array}
 \eeq
 
  Suppose $0<\gamma qr<n$ and $0<\delta \left({p\over p-1}\right)r<n$. 
 We have
\bel{Decay Case3 Est1}
\begin{array}{lr}\ds
\prod_{i=1}^2 \vol\{\Q_i\}^{{\zeta\over n}-{1\over p}+{1\over q}}\left\{{1\over\vol\{\Q_1\}\vol\{\Q_2\}}\iint_{\Q_1\times\Q_2} \left({1\over|u|+|v|}\right)^{\gamma qr}dudv\right\}^{1\over qr}
\\ \ds~~~~~~~~~~~~~~~~~~~~~~~~~~~~
\left\{{1\over\vol\{\Q_1\}\vol\{\Q_2\}}\iint_{\Q_1\times\Q_2} \left({1\over |u|+|v|}\right)^{\delta{pr\over p-1}}dudv\right\}^{p-1\over pr}
\\\\ \ds
=~ \lambda^{\zeta-{n\over p}+{n\over q}}\left({1\over \lambda}\right)^{n\over qr}  \left({1\over \lambda}\right)^{n\left({p-1\over pr}\right)}
\\ \ds~~~~
\left\{\int_{\Q_1}\left\{\int_{\Q_2}  \left({1\over |u|+|v|}\right)^{\gamma qr} dv\right\} du\right\}^{1\over qr}
\left\{\int_{\Q_1}\left\{\int_{ \Q_2}  \left({1\over |u|+|v|}\right)^{\delta{pr\over p-1}} dv\right\}du\right\}^{p-1\over pr}

\\\\ \ds
\leq~ \lambda^{\zeta-{n\over p}+{n\over q}}\left({1\over \lambda}\right)^{n\over qr}  \left({1\over \lambda}\right)^{n\left({p-1\over pr}\right)} \lambda^{n\over qr} \lambda^{n\left({p-1\over pr}\right)}
\left\{\int_{\Q_2}  \left({1\over |v|}\right)^{\gamma qr} dv\right\}^{1\over qr}
\left\{\int_{ \Q_2}  \left({1\over |v|}\right)^{\delta{pr\over p-1}} dv\right\}^{p-1\over pr}
\\\\ \ds
\leq~ \lambda^{\zeta-{n\over p}+{n\over q}}\left({1\over \lambda}\right)^{n\over qr}  \left({1\over \lambda}\right)^{n\left({p-1\over pr}\right)} \lambda^{n\over qr} \lambda^{n\left({p-1\over pr}\right)}
\left\{\int_{\Q_2^*}  \left({1\over |v|}\right)^{\gamma qr} dv\right\}^{1\over qr}
\left\{\int_{ \Q_2^*}  \left({1\over |v|}\right)^{\delta{pr\over p-1}} dv\right\}^{p-1\over pr}
\\ \ds~~~~~~~~~~~~~~~~~~~~~~~~~~~~~~~~~~~~~~~~~~~~~~~~~~~~~~~~~~~~~~~~~~~~~~~~~~~~~~~~~~~~~~~~~~~~~~~
 \hbox{\small{by (\ref{Int Compara gamma})-(\ref{Int Compara delta})}}
\\\\ \ds
\leq~\C_{p~q~r~\gamma~\delta} ~\lambda^{\zeta-{n\over p}+{n\over q}}\left({1\over \lambda}\right)^{n\over qr}  \left({1\over \lambda}\right)^{n\left({p-1\over pr}\right)}
\lambda^{n\over qr}\lambda^{n\left({p-1\over pr}\right)}

\\\\ \ds
=~\C_{p~q~r~\gamma~\delta}~ \lambda^{\zeta-{n\over p}+{n\over q}}
\\\\ \ds
=~\C_{p~q~r~\gamma~\delta}~\lambda^{\gamma+\delta\over 2} \qquad\hbox{( \small{${\zeta\over n}={1\over p}-{1\over q}+{\gamma+\delta\over 2n}$} )}
\\\\ \ds
=~\C_{p~q~r~\gamma~\delta}~\lambda^\ve,\qquad \hbox{\small{$\ve={\gamma+\delta\over 2}>0$}}.
\end{array}
\eeq

 Suppose $n<\gamma qr<2n$ and $0<\delta \left({p\over p-1}\right)r<n$. 
 We have
\bel{Decay Case3 Est2}
\begin{array}{lr}\ds
\prod_{i=1}^2 \vol\{\Q_i\}^{{\zeta\over n}-{1\over p}+{1\over q}}\left\{{1\over\vol\{\Q_1\}\vol\{\Q_2\}}\iint_{\Q_1\times\Q_2} \left({1\over|u|+|v|}\right)^{\gamma qr}dudv\right\}^{1\over qr}
\\ \ds~~~~~~~~~~~~~~~~~~~~~~~~~~~~
\left\{{1\over\vol\{\Q_1\}\vol\{\Q_2\}}\iint_{\Q_1\times\Q_2} \left({1\over |u|+|v|}\right)^{\delta{pr\over p-1}}dudv\right\}^{p-1\over pr}
\\\\ \ds
=~ \lambda^{\zeta-{n\over p}+{n\over q}}\left({1\over \lambda}\right)^{n\over qr}  \left({1\over \lambda}\right)^{n\left({p-1\over pr}\right)}
\\ \ds~~~~
\left\{\int_{\Q_1}\left\{\int_{\Q_2}  \left({1\over |u|+|v|}\right)^{\gamma qr} dv\right\} du\right\}^{1\over qr}
\left\{\int_{\Q_1}\left\{\int_{ \Q_2}  \left({1\over |u|+|v|}\right)^{\delta{pr\over p-1}} dv\right\}du\right\}^{p-1\over pr}
\\\\ \ds
\leq~ \lambda^{\zeta-{n\over p}+{n\over q}}\left({1\over \lambda}\right)^{n\over qr}  \left({1\over \lambda}\right)^{n\left({p-1\over pr}\right)}
\left\{\int_{\Q_1}\left\{\int_{\R^n}  \left({1\over |u|+|v|}\right)^{\gamma qr} dv\right\} du\right\}^{1\over qr}
\left\{\int_{\Q_1}\left\{\int_{ \Q_2}  \left({1\over |v|}\right)^{\delta{pr\over p-1}} dv\right\}du\right\}^{p-1\over pr}
\\\\ \ds
\approx~ \lambda^{\zeta-{n\over p}+{n\over q}}\left({1\over \lambda}\right)^{n\over qr}  \left({1\over \lambda}\right)^{n\left({p-1\over pr}\right)} \lambda^{n\left({p-1\over pr}\right)}
\\ \ds~~~
\left\{\int_{\Q_1}\left\{\idotsint_{\R^n}  \left({1\over |u|+|v_1|+\cdots+|v_n|}\right)^{\gamma qr} dv_1\cdots dv_n\right\} du\right\}^{1\over qr}
\left\{\int_{ \Q_2}  \left({1\over |v|}\right)^{\delta{pr\over p-1}} dv\right\}^{p-1\over pr}
\\\\ \ds
\leq~\C_{p~q~r~\gamma~\delta}~ \lambda^{\zeta-{n\over p}+{n\over q}}\left({1\over \lambda}\right)^{n\over qr}  \left({1\over \lambda}\right)^{n\left({p-1\over pr}\right)} \lambda^{n\left({p-1\over pr}\right)}
\left\{\int_{\Q_1} \left({1\over |u|}\right)^{\gamma qr-n}  du\right\}^{1\over qr}
\left\{\int_{ \Q_2}  \left({1\over |v|}\right)^{\delta{pr\over p-1}} dv\right\}^{p-1\over pr}
\\\\ \ds
\leq~\C_{p~q~r~\gamma~\delta} ~\lambda^{\zeta-{n\over p}+{n\over q}}\left({1\over \lambda}\right)^{n\over qr}  \left({1\over \lambda}\right)^{n\left({p-1\over pr}\right)} \lambda^{n\left({p-1\over pr}\right)}
\left\{\int_{\Q_1^*} \left({1\over |u|}\right)^{\gamma qr-n}  du\right\}^{1\over qr}
\left\{\int_{ \Q_2^*}  \left({1\over |v|}\right)^{\delta{pr\over p-1}} dv\right\}^{p-1\over pr}
\\ \ds~~~~~~~~~~~~~~~~~~~~~~~~~~~~~~~~~~~~~~~~~~~~~~~~~~~~~~~~~~~~~~~~~~~~~~~~~~~~~~~~~~~~~~~~~~~~~~~~~~~~~~~~~~~~~
 \hbox{\small{by (\ref{Int Compara gamma})-(\ref{Int Compara delta})}}

\\ \ds
\leq~\C_{p~q~r~\gamma~\delta}~ \lambda^{\zeta-{n\over p}+{n\over q}}\left({1\over \lambda}\right)^{n\over qr}  \left({1\over \lambda}\right)^{n\left({p-1\over pr}\right)}    \lambda^{n\left({p-1\over pr}\right)}       \lambda^{{2n\over qr}-\gamma}
\\\\ \ds
=~\C_{p~q~r~\gamma~\delta}~\lambda^{\zeta-{n\over p}+{n\over q}}\lambda^{{n\over qr}-\gamma}
\\\\ \ds
~=~\C_{p~q~r~\gamma~\delta}~\lambda^{\gamma+\delta\over 2}\lambda^{{n\over qr}-\gamma} \qquad\hbox{( \small{${\zeta\over n}={1\over p}-{1\over q}+{\gamma+\delta\over 2n}$} )}
\\\\ \ds
=~\C_{p~q~r~\gamma~\delta}~\lambda^{{n\over qr}-{\gamma\over 2}}\lambda^{\delta\over2}
\\\\ \ds
~=~\C_{p~q~r~\gamma~\delta}~\lambda^\ve,\qquad\hbox{\small{$\ve={n\over qr}-{\gamma\over 2}+{\delta\over2}>0$}}.
\end{array}
\eeq  
For $0<\gamma qr<n$ and $n<\delta \left({p\over p-1}\right)r<2n$, an analogue estimate to (\ref{Decay Case3 Est2}) shows the same result with $\ve={\gamma\over2}+n\left({p-1\over pr}\right)-{\delta\over 2}>0$.

Suppose $n<\gamma qr<2n$ and $n<\delta \left({p\over p-1}\right)r<2n$. 
 We have
\bel{Decay Case3 Est3}
\begin{array}{lr}\ds
\prod_{i=1}^2 \vol\{\Q_i\}^{{\zeta\over n}-{1\over p}+{1\over q}}\left\{{1\over\vol\{\Q_1\}\vol\{\Q_2\}}\iint_{\Q_1\times\Q_2} \left({1\over|u|+|v|}\right)^{\gamma qr}dudv\right\}^{1\over qr}
\\ \ds~~~~~~~~~~~~~~~~~~~~~~~~~~~~
\left\{{1\over\vol\{\Q_1\}\vol\{\Q_2\}}\iint_{\Q_1\times\Q_2} \left({1\over |u|+|v|}\right)^{\delta{pr\over p-1}}dudv\right\}^{p-1\over pr}
\\\\ \ds
=~ \lambda^{\zeta-{n\over p}+{n\over q}}\left({1\over \lambda}\right)^{n\over qr}  \left({1\over \lambda}\right)^{n\left({p-1\over pr}\right)}
\\ \ds~~~~
\left\{\int_{\Q_1}\left\{\int_{\Q_2}  \left({1\over |u|+|v|}\right)^{\gamma qr} dv\right\} du\right\}^{1\over qr}
\left\{\int_{\Q_1}\left\{\int_{ \Q_2}  \left({1\over |u|+|v|}\right)^{\delta{pr\over p-1}} dv\right\}du\right\}^{p-1\over pr}
\\\\ \ds
\leq~ \lambda^{\zeta-{n\over p}+{n\over q}}\left({1\over \lambda}\right)^{n\over qr}  \left({1\over \lambda}\right)^{n\left({p-1\over pr}\right)}
\\ \ds~~~~
\left\{\int_{\Q_1}\left\{\int_{\R^n}  \left({1\over |u|+|v|}\right)^{\gamma qr} dv\right\} du\right\}^{1\over qr}
\left\{\int_{\Q_1}\left\{\int_{ \R^n}  \left({1\over |u|+|v|}\right)^{\delta{pr\over p-1}} dv\right\}du\right\}^{p-1\over pr}
\\\\ \ds
\approx~ \lambda^{\zeta-{n\over p}+{n\over q}}\left({1\over \lambda}\right)^{n\over qr}  \left({1\over \lambda}\right)^{n\left({p-1\over pr}\right)}
\\ \ds~~~~~~~
\left\{\int_{\Q_1}\left\{\idotsint_{\R^n}  \left({1\over |u|+|v_1|+\cdots+|v_n|}\right)^{\gamma qr} dv_1\cdots dv_n\right\} du\right\}^{1\over qr}
\\ \ds~~~~~~~
\left\{\int_{\Q_1}\left\{\idotsint_{ \R^n}  \left({1\over |u|+|v_1|+\cdots+|v_n|}\right)^{\delta{pr\over p-1}} dv_1\cdots dv_n\right\}du\right\}^{p-1\over pr}
\\\\ \ds
\leq~\C_{p~q~r~\gamma~\delta}~ \lambda^{\zeta-{n\over p}+{n\over q}}\left({1\over \lambda}\right)^{n\over qr}  \left({1\over \lambda}\right)^{n\left({p-1\over pr}\right)}
\left\{\int_{\Q_1} \left({1\over |u|}\right)^{\gamma qr-n}  du\right\}^{1\over qr}
\left\{\int_{ \Q_1}  \left({1\over |u|}\right)^{\delta{pr\over p-1}-n} du\right\}^{p-1\over pr}
\\\\ \ds
\leq~\C_{p~q~r~\gamma~\delta} ~\lambda^{\zeta-{n\over p}+{n\over q}}\left({1\over \lambda}\right)^{n\over qr}  \left({1\over \lambda}\right)^{n\left({p-1\over pr}\right)}
\left\{\int_{\Q_1^*} \left({1\over |u|}\right)^{\gamma qr-n}  du\right\}^{1\over qr}
\left\{\int_{ \Q_1^*}  \left({1\over |u|}\right)^{\delta{pr\over p-1}-n} du\right\}^{p-1\over pr}
\\ \ds~~~~~~~~~~~~~~~~~~~~~~~~~~~~~~~~~~~~~~~~~~~~~~~~~~~~~~~~~~~~~~~~~~~~~~~~~~~~~~~~~~~~~~~~~~~~~~~~~~~~~~
 \hbox{\small{by (\ref{Int Compara gamma})-(\ref{Int Compara delta})}}

\\ \ds
\leq~\C_{p~q~r~\gamma~\delta}~ \lambda^{\zeta-{n\over p}+{n\over q}}\left({1\over \lambda}\right)^{n\over qr}  \left({1\over \lambda}\right)^{n\left({p-1\over pr}\right)}\lambda^{{2n\over qr}-\gamma}\lambda^{2n\big({p-1\over pr}\big)-\delta}
\\\\ \ds
=~\C_{p~q~r~\gamma~\delta}~\lambda^{\zeta-{n\over p}+{n\over q}}\lambda^{{n\over qr}-\gamma}\lambda^{n\big({p-1\over pr}\big)-\delta} 
\\\\ \ds
=~\C_{p~q~r~\gamma~\delta}~\lambda^{\gamma+\delta\over 2}\lambda^{{n\over qr}-\gamma}\lambda^{n\big({p-1\over pr}\big)-\delta} \qquad\hbox{( \small{${\zeta\over n}={1\over p}-{1\over q}+{\gamma+\delta\over 2n}$} )}
\\\\ \ds
=~\C_{p~q~r~\gamma~\delta}~\lambda^{{n\over qr}-{\gamma\over 2}}\lambda^{n\big({p-1\over pr}\big)-{\delta\over2}}
\\\\ \ds
=~\C_{p~q~r~\gamma~\delta}~\lambda^\ve,\qquad\hbox{\small{$\ve={n\over qr}-{\gamma\over 2}+n\Big({p-1\over pr}\Big)-{\delta\over2}>0$}}.
\end{array}
\eeq

{\small School of Mathematical Sciences, Zhejiang University}\\
{\small sunchuhan@zju.edu.cn}

{\small  Department of Mathematics, Westlake University}\\
{\small wangzipeng@westlake.edu.cn}

\end{document}